\documentclass[12pt,a4paper,reqno]{amsart}

\usepackage[active]{srcltx}

\usepackage{amsmath}
\usepackage{amsfonts}
\usepackage{latexsym}
\usepackage{graphicx}
\usepackage{psfrag}
\usepackage{enumerate}
\usepackage{color}
\usepackage{comment}

\newcommand{\N}{\ensuremath{\mathbb{N}}}
\newcommand{\Z}{\ensuremath{\mathbb{Z}}}
\newcommand{\Q}{\ensuremath{\mathbb{Q}}}

\newcommand{\C}{\ensuremath{\mathbb{C}}}

\def\ds{\displaystyle}
\newcommand{\rme}{\mathrm{e}}

\newcommand{\eps}{\varepsilon}

\input xy
\xyoption{all}
\def\sympow{{\setbox0\hbox{$\bigcirc$}\setbox1\hbox to\wd0{\hss$s$\hss}%
\wd1 0pt\box1\box0}}%symmetric power

\newtheorem{theorem}{Theorem}[section]
\newtheorem{lemma}[theorem]{Lemma}
\newtheorem{proposition}[theorem]{Proposition}
\newtheorem{cor}[theorem]{Corollary}
\newtheorem{remark}[theorem]{Remark}

\newtheorem{corollary}[theorem]{Corollary}

\def\ds{\displaystyle}
\def\rme{\rm e}
\def\proof{\noindent {\sf Proof.}  }
\def\qed{\hfill $\Box$ \\ \bigskip}
\def\K{\mathcal{K}}
\def\L{\mathcal{L}}
\def\Tt{\mathcal{T}}
\def\Rr{\mathcal{R}}
\def\Ss{\mathcal{S}}
\def\Bb{\mathcal{B}}

\date{}

\title[Planar Polynomial Fields and Picard-Vessiot Theory]{On the integrability of polynomial fields in the plane by means of Picard-Vessiot Theory}

\author[P.~Acosta-Hum\'anez]{Primitivo B. Acosta-Hum\'anez}
\address{Department of Mathematics, Universidad del Norte, Colombia}
\email{pacostahumanez@uninorte.edu.co}

\author[J.~T. L\'azaro]{J.~Tom\'as L\'azaro}
\address{Departament de Matem\`atica Aplicada I, Universitat Polit\`ecnica de Catalunya, Spain.}
\email{jose.tomas.lazaro@upc.edu}

\author[J.J.~Morales-Ruiz]{Juan J. Morales-Ruiz}
\address{Technical University of Madrid, Spain.}
\email{juan.morales-ruiz@upm.es}

\author[Ch.~Pantazi]{Chara Pantazi}
\address{Departament de Matem\`atica Aplicada I, Universitat Polit\`ecnica de Catalunya, Spain.}
\email{chara.pantazi@upc.edu}

\thanks{All the authors are partially supported by the MICIIN/FEDER grant number MTM2009-06973 and by the
 Generalitat de Catalunya grant number 2009SGR859. C.P. is additionally partial supported by the MICIIN/FEDER grant  MTM2008--
03437.}



\begin{document}

\begin{abstract}
 We study the integrability of polynomial vector fields using Galois theory of linear differential equations when the
 associated foliations is reduced to a Riccati type foliation. In particular we obtain integrability results for some families of
  quadratic vector fields, Li\'enard equations and equations related with special functions such as Hypergeometric and Heun
  ones. The Poincar\'e problem for some families is also approached.
\end{abstract}

\maketitle

%\vspace{2cm}

\section*{Introduction}

Given a {\it{polynomial differential
system}}\index{polynomial differential system} in $\C^2$,
\begin{equation}
\frac{dx}{dt}={\dot x}=P(x,y), \qquad \qquad \frac{dy}{dt}={\dot
y}=Q(x,y), \label{1}
\end{equation}
with $P,Q \in \C[x,y]$, we consider
its associated differential vector field
\begin{equation}
X= P(x,y) \frac\partial {\partial x}+Q(x,y)\frac \partial {\partial y},
\label{CAMPO}
\end{equation}
whose integral curves are intimately related to the solutions of system~\eqref{1}. These solutions, taken as curves on the plane and leaving for a while its time-dependence, constitute its so-called \emph{foliation} and satisfy the first order differential equation
\begin{equation}
y'  = \frac{dy}{dx} = \frac {Q(x,y)}{P(x,y)}. \label{first}
\end{equation}
This expression~\eqref{first} is often written as a Pfaff equation
\begin{equation}
\Omega =0, \label{form}
\end{equation}
where $\Omega=Q(x,y)dx-P(x,y)dy$ is the corresponding differential $1$-form.
The connection between integral curves of the vector field $X$ and solutions of $\Omega=0$ is clear:
\begin{itemize}
\item geometrically, it is given by $L_X\Omega=0$, which means that the vector field $X$ is tangent to the leaves of the foliation
(the orbits) defined by~\eqref{form};
\item dynamically,
the general solution of equation (\ref{form}), $H(x,y)=\mathit{ctt}$, is given by a {\emph {first integral}} $H$ of the
original vector field $X$, that is, a non-constant scalar function which remains constant along any of its solutions $(x(t),y(t))$.
Since $L_X \Omega =0$, this is equivalent to say that $X(H)=0$ and to the existence of a suitable scalar function $f$ such that
$\Omega=fdH$.
\end{itemize}
From a geometrical point of view, we focus our attention on \emph{invariant algebraic curves}, that is, polynomial integral curves of the vector field $X$. This is the natural framework where Darboux Theory can be applied.
Dynamically speaking, we restrict ourselves to first integrals $H(x,y)$ which are \emph{Liouvillian}, i.e., written
as a combination of algebraic functions, quadratures and exponential of quadratures in ${\C}(x,y)$, the field of rational functions on $x,y$.
As we will see later, Galois Theory provides very useful and powerful tools to approach it.

Two classical problems remain still open for complex polynomial fields:
\begin{itemize}
\item[$(i)$] Concerning the existence of invariant algebraic curves of system~\eqref{1} or, equivalently, of algebraic solutions of the foliation equation~\eqref{first}.
\item[$(ii)$] About the existence of Liouvillian first integrals for systems~\eqref{1} (or, in other words, to determine when
the general solution of equation~\eqref{first} is Liouvillian).
\end{itemize}
For general polynomial vector fields, problems $(i)$ and $(ii)$ are very difficult and we are still far from obtaining an effective method
to decide whether a given arbitrary polynomial field has or not an invariant curve or admits a Liouvillian first integral. In fact, problem $(i)$ is connected with the (also unsolved) classical \emph{Poincar\'e problem}, which seeks for a
bound of the degree of the invariant algebraic curves as a function of the degree of the vector field (or of the associated foliation defined by (\ref{first})).
It is known that Darboux Theory and adjacent results, as the ones due to Prelle--Singer and Singer~\cite{PS,Si}, provide connections between this two problems.

In this work we are concerned with the study of the (Darboux, Galois) integrability of some families of equations of type~\eqref{form} inside the complex analytic category, that is, when the original vector field $X$ defining $\Omega$ is complex polynomial or can be reduced to it. Precisely, we restrict ourselves to those systems which can be reduced to a Riccati type equation
\begin{equation}
v^\prime=a_0(x)+a_1(x)v+a_3(x)v^2, \label{Ric}
\end{equation}
$a_1,a_2$ and $a_3$ being rational functions with complex coefficients.
For Riccati equations there is a very nice theory of integrability in the context of the Galois theory of
its associated second order linear differential equation. This is, in our opinion, a natural framework where several known results
concerning integrability of Riccati  equations~\eqref{Ric} should be considered. Kovacic  provided in 1986 (see \cite{ko}) an
effective algorithm which allows to decide whether an equation\eqref{Ric} has got an algebraic solution or not. And, additionally,
a theorem of Liouville (see~\cite{kol}) proved that
the existence of an algebraic solution is exactly the definition of the integrability for~\eqref{Ric} in the context of
the Galois theory for linear differential equations.
Thus, for foliations of type~\eqref{Ric} problems $(i)$ and $(ii)$ are equivalent
and Kovacic algorithm becomes an extremely powerful tool to approach them.

 In some sense, this work can be considered as a very particular case of the Malgrange approach to the
Galois theory of codimension-$1$ foliations \cite{Malgrange1, Malgrange2,Casale}, that is, for Riccati codimension-$1$ foliations on the complex plane.
Our target is not to obtain general theoretical classification results, but several effective criteria of integrability for such foliations.
We will provide integrability criteria for some families of polynomial quadratic vector fields 
and some  Li\'enard equations involving special functions, allowing this to recover previous results established by several authors.
For instance, we will solve completely the integrability problem for the family of Li\'enard type equations of the form (see \cite{PZ})
\[
yy' = (a(2m+k)x^{2k}+b(2m-k)x^{m-k-1})y-(a^2mx^{4k}+cx^{2k}+b^2m)x^{2m-2k-1},
\]
with $a,b,c,m,k$ complex parameters. It is known that when this
equation comes from a polynomial vector field the constants $m$ and $k$ must
be rational numbers.

We also approach the Poincar\'e problem for some particular families
of systems (see Proposition \ref{orthogonalth} and Theorem \ref{teoremita1}).
\medskip

The paper is structured as follows. To make it as self-contained as possible, we introduce in Section~\ref{definitions} a basic background about Galois theory of linear
differential equations and Darboux theory of integrability of polynomial vector fields.
In Section \ref{someremarks} we remind some useful properties concerning Riccati equations and Section~\ref{applications} is devoted to applications. For completeness we include two Appendixes about
Kovacic algorithm and some special functions.

%%%%%%%%%%%%%%%%%%%%%%%%%%%%%%%%%%%%%%%%%%%%%%%%%%%%%

\section{Two notions of integrability for planar polynomial vector fields}\label{definitions}

\subsection{Darboux theory of Integrability}

We give a very brief overview of Darboux's
integrability ideas \cite{Da}, his terminology and some essential
results.

Let us consider a vector field \eqref{CAMPO} and an irreducible polynomial
$f\in\C[x,y]$.  The curve $f=0$  is called  an \emph{invariant
algebraic curve} of  vector field \eqref{CAMPO} if it satisfies
\[
\dot{f}|_{f=0} = 0.
\]
This condition is equivalent to the existence of a polynomial $K\in \C[x,y]$, called \emph{cofactor}, such that
\[
 X(f(x,y))=P(x,y) \frac{\partial f}{\partial x} + Q(x,y) \frac{\partial f}{\partial y} = K(x,y) f(x,y),
\]
or, equivalently, that
\begin{equation}
\dfrac{X(f)}{f}=X(\log(f))=K.
\label{cofactor}
\end{equation}
From this expression it follows that the curve $f=0$ is formed by
leaves and critical points of the vector field $X=(P, Q)$ defined by
\eqref{CAMPO}.  If the polynomial system~\eqref{1}
has degree $d$, that is $d=\max\left\{ \textrm{deg} P, \textrm{deg}
Q \right\}$, then we have that $\textrm{deg} K \leq d-1$,
independently of the degree of the curve $f(x,y)=0$. From definition
\eqref{cofactor} it follows that if the cofactor $K$ vanishes identically then the polynomial $f$ is a
first integral of the vector field $X$. In terms of the associated foliation, this invariant curve $f=0$ is a particular solution of
$y'=Q/P$ and $Qdx-Pdy=0$.

%to check that $H(x,y)$ is a first integral of the vector field $X$ if and only if $H(x,y)=0$
%is an invariant  curve of $X$ with cofactor $K\equiv 0$.
An analytic $\C$-valued non-constant function $\mu$ is called an
\emph{integrating factor} of system~\eqref{1} if the expression
\[
\dfrac{X(\mu)}{\mu}=X(\log(\mu))=- \nabla \cdot X,
\]
holds,
where $\nabla\cdot X = (\partial P/\partial x) + (\partial
Q/\partial y)$ is the \emph{divergence} of the vector field
$X=(P,Q)$. In case that the domain of definition of $X$ is simply
connected, from the integrating factor $\mu$ it follows that
\[
H(x,y)=\int \mu(x,y) P(x,y)\, dy + \varphi(x),
\]
is a first integral of $X$, provided that $\partial H/\partial y = - \mu Q$.

To ensure the existence of a first integral for a system \eqref{1}
is, in general, a very difficult problem. In~\cite{Da}, Darboux
introduced a method to detect and  construct first integrals using
invariant algebraic curves. Namely, he proved  that any planar
polynomial differential system of degree $d$ having, at least,
$d(d+1)/2$ invariant algebraic curves, admits a first integral or an
integrating factor which can be obtained from them (see also
Jouanolou~\cite{Jo} for a study in  a general context for
codimension-$1$ foliations). Furthermore, Darboux's original ideas
have been improved by taking into account the multiplicity of the
invariant algebraic curves (see \cite{CLP} for more details).
Related to them  some other invariant objects have been introduced
(see \cite{Colin}). They are the so-called \emph{exponential
factors}: given $h,g \in \C[x, y]$ relatively prime, the function
$F={\exp}\left( g/h \right)$ is called an {\emph{exponential
factor}} of the polynomial system (\ref{1}) if there exists a
polynomial $\widetilde{K}\in \C[x,y]$ (also called \emph{cofactor})
that satisfies the equation
\begin{equation}
\dfrac{X(F)}{F}%=\dfrac{X\left(\exp \left( \frac{g}{h} \right)\right)}{\exp \left(
%\frac{g}{h}\right)}
=X\left( \frac{g}{h}\right)=\widetilde{K}. \label{exp}
\end{equation}
It is known that if $h$ is not a constant polynomial then $h=0$ is
an invariant algebraic curve of~\eqref{1} of cofactor $K_h$ satisfying that $X(g)=gK_h + h \widetilde{K}$.

The following theorem (coming originally from Darboux) shows how the
construction of first integrals and integrating factors of
\eqref{CAMPO} can be carried  out from its invariant algebraic curves.

\begin{theorem}
\label{fr} Let consider a planar polynomial system \eqref{1} of
degree $m$, having
\begin{itemize}
 \item $p$ invariant algebraic curves $f_i=0$ with cofactors $K_i$, for $i=1,\ldots,p$ and
 \item $q$ exponential factors $F_j=\exp(g_j/h_j)$ with cofactors $\widetilde{K}_j$, $j=1,\ldots,q$.
\end{itemize}
Then the following assertions hold:
\begin{itemize}
 \item[$(a)$] There exist constants $\lambda_i,\widetilde{\lambda}_j \in \C$ not all vanishing such that
\[
\sum_{i=1}^{p} \lambda_i K_i +  \sum_{j=1}^{q} \widetilde{\lambda}_j
\widetilde{K}_j = 0,
\]
if and only if the multivalued function
\begin{equation}
f_1^{\lambda_1} \ldots f_p^{\lambda_p} F_1^{\widetilde{\lambda}_1}
\ldots F_q^{\widetilde{\lambda}_q}, \label{int}
\end{equation}
is a (Darboux) first integral of system~\eqref{1}.

\item[$(b)$] There exist constants $\lambda_i,\widetilde{\lambda}_j \in \C$ not all vanishing such that
\[
\sum_{i=1}^{p} \lambda_i K_i +  \sum_{j=1}^{q} \widetilde{\lambda}_j
\widetilde{K}_j + \textrm{div}\, X= 0,
\]
if and only if the function defined by~\eqref{int} is a (Darboux) integrating factor of $X$.
\end{itemize}
\end{theorem}
For more recent versions of Theorem \ref{fr} see \cite{LZ08, LZ09}
and  for some  generalizations see \cite{LlPantazi,BP1}.

Functions of the form \eqref{int} are called \emph{Darboux functions}. We say that the polynomial system \eqref{1} is
\emph{Darboux integrable} if it admits a first integral or an integrating factor which is given by a Darboux function.

\begin{remark}\label{remarksinger}
Prelle and Singer \cite{PS} showed that if  system~\eqref{1} admits
an  \emph{elementary first integral} then it admits an integrating
factor which is the $n$-th root of a rational function (a slightly
improved version of this result can be found in~\cite[Corollary
4]{LP1}). Later, Singer in ~\cite{Si} showed that if  system
\eqref{1} admits a Liouvillian first integral  then
 it has an integrating factor which is given by a Darboux function.
This is an important argument to motivate sentences like {``{\it
Darboux functions capture Liouvillian integrability}''} or {``\it
Liouvillian first integrals are either Darboux first integrals or
integrals coming from a Darboux integrating factor}''.
\end{remark}
Given a polynomial system~\eqref{1} of degree $m$, the computation
of all its invariant algebraic curves becomes a complicated problem
since nothing is known a priori about the maximum degree of these
curves. This makes necessary to impose additional conditions either
on the structure of the system ~\eqref{1} or on the {nature} of such
curves (see for instance,~\cite{PER,Car, CeLi,CLPZ} or references
therein). This difficulty has motivated the study of different types
of \emph{inverse problems} of the Darboux theory of integrability
\cite{Pan, CLPZ,CLPW, CLPW2, CLPW3, LPW}.

We finish this subsection with a remark about a geometrical meaning of the integrating factor, pretty known to people coming from the Sophus Lie mathematical community. It is an established fact in Fluid Dynamics that integrating factors arise as a density in stationary planar regimes: the equation $\nabla\cdot (\mu X)=0$ is the continuity equation for the field of velocities $X$, with density function $\mu=\mu(x,y)$ (considered in the context of the Symplectic geometry and Hamiltonian dynamics).
\begin{remark}\label{remarkfactor}
Let $\mu= \mu(x,y)$ be an integrating factor of the vector field \eqref{CAMPO}
defined in some domain of the plane. Then  the vector field $X$ is a hamiltonian vector field with respect to the symplectic form $\Omega=\mu dx\wedge dy$ (this form degenerates only at the zeros of $\mu$). In fact from $i_X\Omega=dH$, we obtain
$$
-\mu Q=\Omega \left(X,\frac{\partial}{\partial x}\right)=\frac{\partial H}{\partial x}, \qquad
\mu P=\Omega \left(X,\frac{\partial}{\partial y}\right)=\frac{\partial H}{\partial y}.
$$
Hence, the vector field \eqref{CAMPO} can be rewritten
\[
 X=P\frac {\partial}{\partial x}+Q\frac {\partial }{\partial y}= \frac 1\mu \left(\frac{\partial H}{\partial y}-\frac{\partial H}{\partial x}\right),
\]
which is hamiltonian with Hamilton function the first integral $H$. It is straightforward to verify that the symplectic form $\Omega$ is invariant under the action of the flow of $X$, i.e.,
\[
  L_X\Omega=di_X\Omega +i_Xd\Omega=di_X\Omega=ddH=0.
\]
From this point of view, the dynamics of the vector fields in the plane can be formally considered as an Ergodic theory problem:
the existence of an invariant measure, the one defined by the associated integrating factor of the flow.
\end{remark}

\subsection{Picard-Vessiot theory}

Picard-Vessiot theory is the Galois theory of linear differential
equations. We will just remind here some of its main definitions and results but we refer the reader to~\cite{vasi} for a wide theoretical background.

We start recalling some basic notions on algebraic groups and,
afterwards, Picard-Vessiot theory will be introduced.

An algebraic group of matrices $2\times 2$ is a subgroup $G\subset
\mathrm{GL}(2,\mathbb{C})$ defined by means of algebraic equations
in its matrix elements and in the inverse of its determinant. That
is, there exists a set of polynomials $P_i \in \C[x_1,\ldots,x_5]$,
for $i\in I$, such that $A\in\mathrm{GL}(2,\mathbb{C})$ given by
\[
A=\left(\begin{array}{cc} x_{11} & x_{12} \\ x_{21} & x_{22}
\end{array}\right),
\]
belongs to $G$ if and only if  $P_i\left( x_{11},x_{12},x_{21},x_{22},\left(\det A \right)^{-1} \right) = 0$ for all $i\in I$ and where $\det A=x_{11}x_{22}-x_{21}x_{12}$.
It is said that $G$ is an algebraic manifold endowed with a
group structure.

Recall that a group $G$ is called \emph{solvable} if and only if there
exists a chain of normal subgroups
\[
e=G_0\triangleleft G_1 \triangleleft \ldots \triangleleft G_n=G,
\]
satisfying that the quotient $G_i/G_j$ is abelian for all $n\geq i\geq j\geq 0$.

It is well known that any algebraic group $G$ has a unique connected normal
algebraic subgroup $G^0$ of finite index. In particular, the
\emph{identity connected component $G^0$} of $G$ is defined as the largest connected
algebraic subgroup of $G$ containing the identity. In case that $G=G^0$ we say that $G$ is a
\emph{connected group}. Moreover, if $G^0$ is solvable we say that $G$ is \emph{virtually solvable}.

The following result provides the relation between virtual solvability of an algebraic group
and its structure.

\begin{theorem}[Lie-Kolchin]\label{LiKo}
Let $G\subseteq \mathrm{GL}(2,\mathbb{C})$ be a virtually solvable group.
Then, $G^0$ is triangularizable, that is, it is conjugate to a subgroup of upper triangular matrices.
\end{theorem}

Now, we briefly introduce Picard-Vessiot Theory.

First, we say that $\left( \K, \phantom{i}' \ \right)$ - or, simply,
$\K$ - is a \emph{differential field} if $\K$ is a commutative field
of characteristic zero, depending on $x$ and $\phantom{i}'$ is a
derivation on $\mathcal{K}$ (that is, satisfying that $(a+b)'=a'+b'$
and $(a\cdot b)'=a'\cdot b+a \cdot b'$ for all $a,b \in \K$). We
denote by $\mathcal{C}$ the \emph{field of constants of $\K$},
defined as $\mathcal{C}=\left\{ c \in \mathcal{K} \ | \ c'=0
\right\}$.

We will deal with second order linear homogeneous differential
equations, that is, equations of the form
\begin{equation}
\label{soldeq}
y''+b_1y'+b_0y=0,\quad b_1,b_0\in \K,
\end{equation}
and we will be concerned with the algebraic structure of their solutions. Moreover, along this work, we will refer
to the current differential field as the smallest one containing the field of coefficients of this differential equation.

Let us suppose that $y_1, y_2$ is a basis of solutions of equation
\eqref{soldeq}, i.e., $y_1, y_2$ are linearly independent over $\K$
and every solution is a linear combination over $\mathcal{C}$ of
these two. Let $\L= \K\langle y_1, y_2 \rangle =\K(y_1, y_2, y_1',
y'_2)$ be the differential extension of $\K$ such that $\mathcal{C}$
is the field of constants for $\K$ and $\L$. In this terms, we say
that $\L$, the smallest differential field containing $\K$ and
$\{y_{1},y_{2}\}$, is the \textit{Picard-Vessiot extension} of $\K$
for the differential equation~\eqref{soldeq}.

The group of all the differential automorphisms of $\L$ over $\K$
that commute with the derivation $\phantom{i}'$  is called the
\emph{Galois group} of $\L$ over $\K$ and is denoted by ${\rm
Gal}(\L/\K)$. This means, in particular, that for any $\sigma\in
\mathrm{Gal}(\L/\K)$, $\sigma(a')=(\sigma(a))'$ for all $a\in \L$
and that $\sigma(a)=a$ for all $a\in \K$. Thus, if $\{y_1,y_2\}$ is
a fundamental system of solutions of~\eqref{soldeq} and $\sigma \in
\mathrm{Gal}(\L/\K)$ then $\{\sigma y_1, \sigma y_2\}$ is also a
fundamental system. This implies the existence of a non-singular
constant matrix
\[
A_\sigma=
\begin{pmatrix}
a & b\\
c & d
\end{pmatrix}
\in \mathrm{GL}(2,\mathbb{C}),
\]
such that
\[
\sigma
\begin{pmatrix}
y_{1}&
y_{2}
\end{pmatrix}
=
\begin{pmatrix}
\sigma (y_{1})&
\sigma (y_{2})
\end{pmatrix}
=\begin{pmatrix} y_{1}& y_{2}
\end{pmatrix}A_\sigma.
\]
This fact can be extended in a natural way to a system
\[
\sigma
\begin{pmatrix}
y_{1}&y_2\\
y'_1&y'_{2}
\end{pmatrix}
=
\begin{pmatrix}
\sigma (y_{1})&\sigma (y_2)\\
\sigma (y'_1)&\sigma (y'_{2})
\end{pmatrix}
=\begin{pmatrix} y_{1}& y_{2}\\y'_1&y'_2
\end{pmatrix}A_\sigma,
\]
which leads to a faithful representation $\mathrm{Gal}(\L/\K)\to
\mathrm{GL}(2,\mathbb{C})$ and makes possible to consider
$\mathrm{Gal}(\L/\K)$ as a subgroup of $\mathrm{GL}(2,\mathbb{C})$
depending (up to conjugacy) on the choice of the fundamental system $\{y_1,y_2\}$.

One of the fundamental results of the Picard-Vessiot Theory is the
following theorem (see~\cite{ka,kol}).

\begin{theorem}  The Galois group $\mathrm{Gal}(\L/\K)$ is an
algebraic subgroup of $\mathrm{GL}(2,\mathbb{C})$.
\end{theorem}

We say that equation~\eqref{soldeq} is \emph{integrable} if the Picard-Vessiot extension
$\L\supset \K$ is obtained as a tower of differential fields
$\K=\L_0\subset \L_1\subset\cdots\subset \L_m=\L$ such that
$\L_i=\L_{i-1}(\eta)$ for $i=1,\ldots,m$, where either
\begin{itemize}
\item[$(i)$] $\eta$ is {\emph{algebraic}} over $\L_{i-1}$, that is $\eta$ satisfies a
polynomial equation with coefficients in $\L_{i-1}$.
\item[$(ii)$] $\eta$ is {\emph{primitive}} over $\L_{i-1}$, that is $\eta' \in \L_{i-1}$.
\item[$(iii)$] $\eta$ is {\emph{exponential}} over $\L_{i-1}$, that is $\eta' /\eta \in \L_{i-1}$.
\end{itemize}

Usually in terms of  Differential Algebra's terminology we say that
equation~\eqref{soldeq} is integrable if the corresponding
Picard-Vessiot extension is \emph{Liouvillian}. Moreover, the
following theorem  holds.

\begin{theorem}[Kolchin]
\label{LK}
Equation~\eqref{soldeq} is integrable if
and only if $\mathrm{Gal}(\L/\K)$ is virtually solvable, that is, its identity component
$(\mathrm{Gal}(\L/\K))^0$ is solvable.
\end{theorem}

For instance, for the case $a=0$ in equation~\eqref{soldeq}, i.e.
$y''+by=0$, it is very well known~\cite{ka,kol,vasi} that
$\mathrm{Gal}(\L/\K)$ is an algebraic subgroup of ${\rm
SL}(2,\mathbb{C})$ (remind that $A\in {\rm SL}(2,\mathbb{C})
\Leftrightarrow A \in {\rm GL}(2,\mathbb{C})$ and $\det A=1$). For a
more detailed study see the Appendix~\ref{kovsection}.

\section{Some remarks about Riccati equation}\label{someremarks}

Ricatti equation is probably one of the most studied equations in Dynamical Systems. Its r\^{o}le in the study of the Darboux and Picard-Vessiot
integrability leads us to devote this section to some of its properties. Even though these results
are known, their proofs have been included for completeness. We divide these properties in two types:
the first one (see Subsection \ref{trns}) concerning transformations leading a general second order differential equation into a Riccati equation
(written in the so-called reduced form). Remind that this becomes the starting point of the celebrated Kovacic algorithm
(see Appendix~\ref{kovsection}). A second type, more Darboux-like, that studies first integrals and integrating factors for a Riccati equation
(Subsection \ref{subint}).

\subsection{Transformations related to Riccati equations} \label{trns}

It is known that any second order differential equation can be led into a general Riccati equation
through a classical logarithmic change of variable (see, for instance,~\cite{PZ, INC}). The following proposition recall it and summarises some other related transformations.

\begin{proposition}\label{prop:almp}
Let $\K$ be a differential field and let consider functions $a_0(x)$,
$a_1(x)$, $a_2(x)$, $r(x)$, $\rho(x)$, $b_0(x)$, $b_1(x)$ belonging to $\K$ that, for simplicity, will be denoted without their explicit dependence on $x$. Consider now the following forms associated to any second order differential equation (ode) and Riccati equation:
\begin{itemize}
\item[$(i)$] Second order ode (in general form):
\begin{equation}
\label{lemma:2ode:initial_form} y''+b_1 y' + b_0 y=0.
\end{equation}

\item[$(ii)$] Second order ode (in reduced form):
\begin{equation}
\label{lemma:2ode:normal_form} \xi''=\rho\xi.
\end{equation}

\item[$(iii)$] Riccati equation (in general form):
\begin{equation}
\label{lemma:riccati:initial_form} v'=a_0 + a_1 v + a_2 v^2
,\quad a_2 \neq 0.
\end{equation}

\item[$(iv)$] Riccati equation (in reduced form):
\begin{equation}
\label{lemma:riccati:normal_form} w'=r - w^2,
\end{equation}
\end{itemize}
Then, there exist transformations $\Tt$, $\Bb$, $\Ss$ and $\Rr$ leading some of these equations into the other ones, as showed in the following diagram:
\[
\xymatrix{ v'=a_0 + a_1 v + a_2 v^2
\ar[r]^-{\Tt} & w'=r - w^2 \\
y''+b_1 y' + b_0 y=0 \ar[r]^-{\Ss} \ar@{<-}[u]^{\Bb} &
\ar@{<-}[u]_{\Rr} \xi''=\rho \xi. }
\]
The new independent variables are defined by means of
\[
\begin{array}{lcl}
\Tt: \ v=- \left( \dfrac{a_2'}{2a_2^2} +
\dfrac{a_1}{2a_2} \right) - \dfrac{1}{a_2}w, &\qquad \qquad &
\Bb: \ {\ds v=-\frac{1}{a_2} \frac{y'}{y} }, \\\\\
\Ss: \ {\ds y=\xi {\rme}^{-{\frac{1}{2}}\int b_1 dx}}, &\qquad \qquad &
\Rr: \ {\ds w=\frac{\xi'}{\xi} },
\end{array}
\]
and the functions $r$, $\rho$, $b_0$ and $b_1$ are given by
\begin{align}
r &=\frac{1}{\beta} \left( a_0 + a_1 \alpha + a_2
\alpha^2 - \alpha' \right), \\ %[1.2ex]
\alpha &=-\left({\frac{a_2'}{2a_2^2}}+{\frac{a_1}{2a_2}}\right),\quad \beta=-\frac1{a_2}, \\ %[1.2ex]
b_1 &= -\left(a_1+{\frac{a_2'}{a_2}}\right),\quad b_0=a_0a_2,  \\ %[1.2ex]
\rho &= r = \frac{b_1^2}{4}+\frac{b_1'}{2} -b_0.
\end{align}

\end{proposition}

\proof The proof is quite standard.

\begin{itemize}
\item[[$\Tt\!\!$]]: Applying the change $v=\alpha +
\beta w$ we get the equation
\[
\alpha'+ \beta' w + \beta w' = a_0 + a_1 \alpha + a_1 \beta w + a_2
\alpha^2+ 2 a_2 \alpha \beta w + a_2 \beta^2 w^2
\]
that, regrouping terms, leads to
\[
w' = \frac{1}{\beta} \left( a_0 + a_1 \alpha + a_2 \alpha^2 -
\alpha' \right) + \left( a_1 + 2a_2 \alpha - \frac{\beta'}{\beta}
\right) w + a_2 \beta w^2.
\]
Since $a_2 \neq 0$ we can take $\beta=-1/a_2$ and, therefore,
$a_2 \beta = -1$. Having this into account, the value of
$\alpha$ satisfying that the coefficient in $w$ vanishes is given by
\[
\alpha=\frac{1}{2a_2} \left( \frac{\beta'}{\beta} - a_1 \right).
\]
The expressions for $\alpha$, $\beta$ and $r$ follow
straightforwardly,
\[
{\ds r=\frac{1}{\beta} \left( a_0 + a_1 \alpha + a_2 \alpha^2 -
\alpha' \right)},\quad \alpha=-\left( \frac{a_2'}{2a_2^2}+\frac{a_1}{2a_2}\right),\quad \beta=-\frac1{a_2}.
\]
Moreover, it is clear that $\alpha$, $\beta$ and $r$ belong to $\K$.

\medskip

\item[[$\Bb\!\!$]]: Imposing $\alpha=0$ and taking $\beta=-1/a_2$ in transformation $\Tt$ we have $v=-w/a_2$ and we obtain the Riccati equation
\[
w'=-a_0a_2+\left(a_1+ \frac{a_2'}{a_2} \right)w-w^2.
\]
Performing now the change of variables $w=(\log y)'$ (or, equivalently, $v=-a_2 y'/y$) we obtain the differential equation
$y''+b_1 y'+b_0 y=0$ with
\[
b_1=-\left(a_1+\frac{a_2'}{a_2}\right), \quad
b_0=a_0a_2.
\]
Obviously, $b_0$ and $b_1$ belong to $\K$.

\medskip

\item[[$\Ss\!\!$]]: The change of variable $y=\mu\xi$, with $\mu=\mu(x)$ and $\xi=\xi(x)$, lead us to
\[
\xi''+\left(2 \frac{\mu'}{\mu} +b_0\right)\xi'+\left( \frac{\mu''}{\mu}+b_0 \frac{\mu'}{\mu}+b_1\right)\xi=0.
\]
In order to obtain the equation $\xi''=\rho \xi$ we need to impose
\[
2\frac{\mu'}{\mu} + b_0=0,\qquad \frac{\mu''}{\mu} + b_0 \frac{\mu'}{\mu} + b_1=-\rho,
\]
which gives rise to
\[
\mu=\rme^{-\frac{1}{2} \int b_0},\qquad
\rho=\frac{b_0^2}{4}+\frac{b'_0}{2}-b_1.
\]
Moreover, it is straightforward to check that $\rho\in \K$.

\medskip

\item[[$\Rr\!\!$]]: This is a particular case of transformation $[\Bb]$ with the particular choice $a_0=r$, $a_1=0$ and $a_2=-1$.
\end{itemize}

\medskip

Finally, composing the transformations provided by $[\Bb]$, $[\Rr]$ and $[\Ss]$:
\[
-a_2v=\frac{y'}{y},\qquad y=\xi e^{-\frac12\int b_0},\qquad
b_0=-\left(a_1+ \frac{a_2'}{a_2} \right)\qquad \frac{\xi'}{\xi}=w,
\]
we recover the result given by $[\Tt]$,
\[
v=-\left(\frac{a_1}{2a_2} + \frac{a_2'}{2a_2^2} \right)-\frac1{a_2}w=\alpha+\beta w,
\]
which implies that, in some sense and taking $\rho=r$, the diagram commutes: $\Ss \circ \Bb=
\Rr \circ \Tt$.
\qed

From this Lemma, it follows that the function $v$ is algebraic over $\K$ if and only if the function $w$ is also algebraic over $\K$. Furthermore, in such case, the degree over $\K$ of both functions $v$ and $w$ is the same.

It is known that a Riccati equation~\eqref{lemma:riccati:initial_form} has an algebraic solution
over $\K$ if and only if the differential equation~\eqref{lemma:2ode:initial_form} is integrable in a Picard-Vessiot sense.
 In this situation we say that the Riccati equation is integrable over $\K$. We notice that
Kovacic algorithm (see Appendix~\ref{kovsection}) starts from an equation in form~\eqref{lemma:2ode:normal_form}.

\subsection{Integrating factor and first integrals for Riccati  vector fields} \label{subint}

We briefly show some relations between the existence of invariant curves of a certain type of vector fields and the integrability, via Kovacic algorithm (see Appendix~\ref{kovsection}), of its associated  Riccati foliation.

%%% Tom\'as: por favor, comprobad si lo que he puesto es correcto o no.

%% Inicio del p\'arrafo %%%%%%%%%%%%%%%%%%%%%%%%%%%

From Singer~\cite{Si} (see Remark \ref{remarksinger}) we know that if a planar polynomial vector field \eqref{CAMPO}
admits a Liouvillian first integral then it has also an integrating factor given by a Darboux function. However, very few results are known
about the relation between the existence of an algebraic invariant curve of a general planar vector field and the Liouvillian integrability of its foliation.

Let us consider a family of planar vector fields of the form
\begin{equation}
X=(p(x)-q(x)w^2)\frac{\partial }{\partial w}+ q(x)\frac{\partial
}{\partial x},
\label{vric}
\end{equation}
with $p(x),q(x) \in \C[x]$ complex polynomials. Introducing an independent variable $t$, usually called \emph{time}, we can  associate to them the following system of differential equations
\[
\begin{array}{rcl}
\dot{w} &=& p(x) - q(x) w^2, \\
\dot{x} &=& q(x),
\end{array}
\]
where we denote by $\dot{\phantom{x}}=d/dt$. Its foliation, governed by the equation
\[
w'=\dfrac{dw}{dx} = \dfrac{p(x) - q(x) w^2}{q(x)} = \dfrac{p(x)}{q(x)} - w^2,
\]
is a Riccati equation given in reduced form $w'=r(x) - w^2$ with $r=p/q \in \C(x)$.
Next lemma shows that the integrability of this ``Riccati foliation'' is closely related to
the existence of an algebraic invariant curve of its vector field~\eqref{vric}. A similar approach for this problem can be
 found in~\cite{GGG0, GGG}.

\begin{lemma}\label{genintfactor}
Let $w_1=w_1(x)$ be a solution of a Riccati equation in reduced form
\[
w'=r(x)-w^2,
\]
with $r(x)=p(x)/q(x)\in \C[x]$. Then the associated vector field \eqref{vric}
has an integrating factor given by
\begin{equation}
\mu_1(w,x)=\dfrac{e^{-2\displaystyle{\int} w_1(x) dx}}{(-w+w_1(x))^2}.
\label{intR1}
\end{equation}
\end{lemma}
\proof It is straightforward to check that if $w_1(x)$ is a solution of
$w'=p/q - w^2$ then it holds $X(f_1)=K_1f_1$ with
$f_1(w,x)=-w+w_1(x)=0$ and $K_1=-q(w+w_1(x))$. In addition,
$F_1(x)=\rme^{-\int\omega_1(x) dx}$ satisfies $X(F_1)=L_1 F_1$ with
$L_1=-q\omega_1.$ Notice that $X$ has divergence $\mbox{div}X=-2qw$, and
additionally it holds $-2K_1+2L_1+\mbox{div}X=0$. Proceeding as in Theorem~\ref{fr}, vector field \eqref{vric} admits the integrating
factor
\[
\mu_1(w,x)=\dfrac{F_1^2}{f_1^2}=\dfrac{e^{-2\int w_1(x) dx}}{(-w+w_1(x))^2},
\]
as it was claimed.
\qed

\begin{remark}
 We  stress that the result in Lemma \ref{genintfactor}
is independent of the nature of the solution $w=w_1(x)$. In general,  the integral $\int\omega_1(x)dx$ is an abelian integral.
\end{remark}

The important fact is that, conversely, Picard-Vessiot theory and in particular, Kovacic algorithm, supply information about
 first integrals and integrating factors of the equation $w'=r(x)-w^2$ from the knowledge of some of its solutions,
 $w_1,w_2,w_3$. Indeed, from the first three cases in Kovacic algorithm~\cite{ko} (the integrable ones)
 one obtains the following types of first integrals (see Weil~\cite{weil} and \.{Z}o{\l}\c{a}dek~\cite{Zol}).

\begin{proposition} \label{propclas}
The following statements hold.
\begin{itemize}
\item[\textsf{Case 1:}] One has two possibilities:
\begin{itemize}

\item[-] If only $w_1\in \C(x)$ then  $X$ admits a first integral of Darboux--Schwarz--Christoffel type.

\item[-] If both $w_1, w_2\in \C(x)$ then  $X$ admits a first integral of Darboux type. In
particular, from Lemma \ref{genintfactor} we can construct two
integrating factors $\mu_1$ and $\mu_2$ and so $\mu_1/\mu_2$ is a first
integral of $X$. Thus, we have
\[
H(w,x)=\dfrac{(-w+w_2(x))}{(-w+w_1(x))}e^{\int[(w_2(x)-w_1(x)) dx]}.
\]
\end{itemize}

\item[\textsf{Case 2:}]  If $w_1$ is a solution of a quadratic
polynomial  then  $X$ admits a first integral  of hyperelliptic type.

\medskip

\item[\textsf{Case 3:}]  If all $w_1,w_2,w_3$ are  algebraic over
$\C(x)$ then $X$ admits a rational first integral.
\end{itemize}
\end{proposition}
The following result characterises the rational integrability of the polynomial
vector fields that we are considering in this work.
\begin{cor}
The Galois group of \eqref{lemma:2ode:normal_form} is finite if and only if
its corresponding planar polynomial vector field has a rational first
integral.
\end{cor}
\proof
Notice that $\mathrm{Gal}(\L/\K)$ is finite if and only if we fall in
case 3 of Kovacic Algorithm or in case 1 of the type
$$
\mathrm{Gal}(\L/\K)=\left\{ \left(
\begin{array}{cc}
c & 0\\
0& c^{-1}
\end{array}
\right), \quad c^n=1 \right\}.
$$
Hence, only remains to study this last case (the cyclic one).
 Let $\xi_1,\xi_2$ be  solutions of $\xi''=\rho \xi.$ Then
there exists $g\in\C(x)$ such that $\xi_1=g^{\frac{1}{n}}$ and
$\xi_2=g^{-\frac{1}{n}}.$ We define $\omega_1=\xi_1'/\xi_1$,
$\omega_2=\xi_2'/\xi_2,$  and we obtain
$$
\omega_1=\frac{1}{n}\frac{g'}{g},  \qquad \omega_2=-\frac{1}{n}\frac{g'}{g}.
$$
Thus, the corresponding vector field \eqref{lemma:2ode:normal_form} admits the first integral
$$
H_1(\omega,x)=\dfrac{-\omega+\omega_2}{-\omega+\omega_1}e^{\int (\omega_2-\omega_1)dx}=\dfrac{-\omega-\dfrac{1}{n}\dfrac{g'}{g}}{-\omega+\dfrac{1}{n}\dfrac{g'}{g}}\ g^{\frac{-2}{n}},
$$
and so also admits the first integral
$$
H(\omega,x)=H_1(\omega,x)^n=\dfrac{1}{g^2}\left(\dfrac{-ng\omega-{g'}}{-ng\omega+{g'}}\right)^n\in\C(\omega,x),
$$
which completes the proof.
\qed

\bigskip

\begin{remark}
Let $P\in\C[x]$ be a polynomial of odd degree. It is known that the planar polynomial vector
field $\dot{x}=1, \  \dot{y}=P(x)+y^2$ (with associated foliation
$y'=P(x)+y^2$) is not integrable, that is, it has no invariant curves,
since it falls in case 4 of Kovacic Algorithm~\cite{ko}.
\end{remark}

\subsection{Riccati foliations}

Let us recall some well-known geometrical properties of Riccati foliations defined by planar polynomials vector fields (see, for instance, \cite{LN}).
Although we are not using these properties along the paper, we include them for completeness.

Let
\begin{equation}
\Omega=q(x)dy-(p_1(x)+p_2(x)y+p_3(x)y^2)dx, \quad  p_i,q\in {\C}[x]\label{Ricf}
\end{equation}
be the $1$-form defining a Riccati foliation on the complex plane. Since the equation $\Omega=0$ is the projective version of the corresponding
second order linear differential equation defined over the vector bundle ${\bf P}^1\times \C^2$ (i.e., the fibre $\C^2$ is projectivized to ${\bf P}^1$, see subsection \ref{trns}), the Riccati equation $\Omega=0$ is defined in a natural way over $(x,y)\in {\bf P}^1\times{\bf P}^1$.

The singular points of $\Omega=0$ are the zeros  of the polynomial $q(x)$ and, possibly,  the point at infinity $x=\infty\in{\bf P}^1$.
Moreover, these singular points are exactly the poles of the coefficients of the associated second order linear differential equation.
 We define $d:=max (deg(p_1), deg(p_2),deg(p_3), deg(q)-2)$. Thus the point $x=\infty $ is a singular point if and only if $deg(q)$ is less than $d+2$.
In  Kovacic  algorithm, which applies to the reduced form of the second order linear differential equation (see Appendix~\ref{kovsection}),
this set of singular points is denoted by $\Gamma$. Therefore it seems natural to call it in the same way also here, that is, $\Gamma=\{x_1,...,x_r\}$.
Thus the Riccati foliation is {\it holomorphic} on $({\bf P}^1-\Gamma)\times{\bf P}^1$, because the singular points of the associated linear differential
equation are exactly the set $\Gamma$, given by the singularities of their coefficients.  We notice that the (``singular") sets $\{x_i\}\times {P}^1$ are invariant by the
foliation and are usually called \emph{invariant fibres} because they are already fibres of the 
 fibration $\pi: {\bf P}^1\times{\bf P}^1\rightarrow {\bf P}^1$,
$(x,y)\mapsto x$. This fibration is  transversal  to the Riccati foliation since fibres $T(x):=\pi^{-1} (x)=\{x\}\times {\bf P}^1$, 
with $x$ non
singular, are {\it global} and transversal to the foliation, i.e., transversal to all the leaves. Over any of these transversals the holonomy group of the foliation is
 defined as
a representation
$$
\pi_1({\bf P}^1-\Gamma, x_0)\rightarrow \mbox{Diff}(T(x_0)),
$$
(where $\mbox{Diff}(T(x_0))$ is the group of diffeomorphisms on the transversal),
given by lifting the loops in the fundamental group to the leaves of the foliation, this is, by solving the Riccati equation with initial conditions and final
points on the transversal $T(x_0)$. As the Riccati equation is the projectivization of a second order linear differential equation, the holonomy group must be the projectivization
of the monodromy group of the linear second order equation acting on the vector space fibre of the meromorphic vector bundle ${\bf P}^1\times {\C}^2$,
$\{x\}\times {\C}^2\approx {\bf C}^2$. By fixing a base of fundamental solutions, this can be considered as the space of solutions of the linear differential equation.
 Hence, as the monodromy group is represented by the linear group $GL(2, \C)$, the holonomy group is represented by the projective linear group $PGL(2, \C)$, the M\"obius
transformations
$$
\pi_1({\bf P}^1-\Gamma, x_0)\rightarrow PGL(2, \C)
.$$
The Riccati foliations are the most well-known class of a family of foliations, the  projective foliations, with holonomy group represented in the projective group. For Riccati foliations, the holonomy group is contained in the Galois group of the foliation either in the Malgrange approach \cite{Malgrange1, Malgrange2,Casale} or in the Lie-Vessiot-Kolchin approach \cite{BM, BS}. In fact,
if the singular points of the associated linear differential equations are singular regular ones, then the Zariski adherence of the holonomy group is the
Galois group of the foliation, because in this case the Zariski adherence of the monodromy group of the linear differential equation is the Galois group of the associated  linear differential equation.

The critical points of the associated vector field, i.e., zeroes of $q(x)$ and of $p_1(x)+p_2(x)y+p_3(x)y^2$ are obviously contained in the invariant fibres. For general Riccati foliations there are two critical points on the invariant fibre for each point in $\Gamma$. For some special Riccati foliations there are no critical points. This is the case, for example, of Riccati foliations given in reduced form by $\Omega=q(x)dy-(p(x)+q(x)y^2)dx$ (with  $p$ and $q$ relatively primes), corresponding to the field (\ref{vric}).

%%%%%%%%%%%%%%%%%%%%%%%%%%%%
%%%%%%%%%%%%%%%%%%%%%%%%%%%%

\section{Applications} \label{applications}
In this section we analyse some examples involving integrability and non-integrability of some families of Riccati planar vector fields or planar vector fields whose foliation can be reduced into a Riccati form.

\subsection{Quadratic polynomials fields}
The study of the integrability of the  quadratic polynomial vector field
$$
\begin{array}{rcl}
\dot{x} &=& a_{20}x^2+a_{11}xy+a_{02}y^2+a_{10}x+a_{01}y+a_{00}, \\
\dot{y} &=& b_{20}x^2+b_{11}xy+b_{02}y^2+b_{10}x+b_{01}y+b_{00},
\end{array}
$$
with $a_{ij},b_{i,j}\in\mathbb{C}$ is, in its general form, a hard problem.
One of its possible approaches is the so-called linear-quadratic case, when one of the two components is a polynomial of degree one.
In~\cite[Prop.3]{LlibreValls} it is proved that its study
around a  finite equilibrium point (the origin)
can be reduced to consider two families of
systems. Using the notation introduced therein, we refer to these families as (S1)-type,
$$
\begin{array}{cr}

\begin{array}{rcl}
\dot{x} &=& x, \\
\dot{y} &=& \eps x + \lambda y + b_{20}x^2 + b_{11}xy + b_{02}y^2,
\end{array}
& \qquad \qquad \mbox{(S1)}
\end{array}
$$
and (S2)-type,
$$
\begin{array}{cr}
\begin{array}{rcl}
\dot{x} &=& y, \\
\dot{y} &=& \eps x + \lambda y + b_{20}x^2 + b_{11}xy + b_{02}y^2.
\end{array}
&
\qquad \qquad \mbox{(S2)}
\end{array}
$$
In \cite{LlibreValls}, the authors prove that the linear-quadratic
systems having a global analytic first integral are those
satisfying:
\begin{itemize}
 \item[$(a_1)$] $b_{02}=\lambda=0$.
 \item[$(b_1)$] $b_{02}=0$ and $\lambda=-p/q \in \Q^{-}$,
\end{itemize}
in the case of (S1)-type systems  and
\begin{itemize}
 \item[$(a_2)$] $b_{20}=b_{02}=\lambda=0$ and $\eps b_{11} \neq 0$.
 \item[$(b_2)$] $b_{20}=b_{11}=\lambda=0$ and $\eps b_{02} \neq 0$.
 \item[$(c_2)$] $b_{11}=\lambda=0$ and $b_{20} \neq 0$,
\end{itemize}
for (S2)-type systems. Furthermore, they also
provide the explicit form of the corresponding first integrals. It
is important to notice that all of them are of Darboux type and, therefore, Liouvillian.

Our aim in this example is to show that these results can be recovered using
arguments coming from the Galois theory of linear differential equations. We
start first with the (S1)-case, whose associated
foliation is given by the Riccati equation:
\begin{equation}
\frac{dy}{dx}=\left( \eps + b_{20}x \right) + \left( \frac{\lambda +
b_{11}x}{x} \right) y + \frac{b_{02}}{x} y^2. \label{S1:Riccati}
\end{equation}
By Lemma \ref{prop:almp} this equation
can be transformed into the reduced form $ w^\prime=r(x)-w^2$, with
\begin{equation}
\label{Whittaker-type} r(x)=\frac{1}{4} - \frac{\kappa}{x} +
\frac{4\mu^2-1}{4x^2}, \quad \kappa =
\frac{1}{\sqrt{b_{11}^2-4b_{20}b_{02}}} \left( b_{02} \eps +
\frac{b_{11}}{2} (1-\lambda) \right), \quad \mu= \frac{\lambda}{2},
\end{equation}
provided $b_{11}^2-4b_{20}b_{02} \neq 0$, and into
the form $\xi''=r(x)\xi$. This equation is a Whittaker equation (see Appendix B) to which one can apply
the Martinet-Ramis Theorem (see theorem~\ref{thmarram}). This Theorem asserts that such Whittaker equation is
integrable if and only if at least one of the following conditions is
verified:
\[
 \pm \kappa \pm \mu \in \frac{1}{2} + \N,
\]
or, equivalently (and more suitable for the expressions derived of
$\kappa$ and $\mu$)
\[
 2 \left( \kappa \pm \mu \right) \in 2\Z + 1.
\]
In our case one has that
\[
2 \left( \kappa \pm \mu \right) = \frac{2b_{20}\eps +
b_{11}(1-\lambda)}{\sqrt{b_{11}^2 - 4b_{20} b_{02}}} \pm \lambda,
\]
so for (S1)-type systems conditions $(a_1)$ and $(b_1)$ read, respectively,
$2(\kappa \pm \mu)= 1 \in 2\Z+1$ and $2(\kappa +
\mu) = (1 + (p/q)) + (-p/q) = 1 \in 2\Z+1$. Therefore Galois Theory
recovers the integrability result asserted
in~\cite[Thm.1]{LlibreValls}.

Let us consider now a (S2)-type system, namely,
\[
\begin{array}{rcl}
\dot{x} &=& y, \\
\dot{y}&=& \eps x + \lambda y + b_{20}x^2 + b_{11}xy + b_{02}y^2,
\end{array}
\]
with foliation given by the differential equation
\begin{equation}
\label{S2:foliation_eq} \frac{dy}{dx} = \left( \lambda + b_{11}x
\right) + \left( \eps x + b_{20} x^2 \right) \frac{1}{y} + b_{02} y.
\end{equation}
This equation falls in one of the following situations:
\begin{itemize}
\item[$(i)$] $\lambda=b_{11}=0$ yields to a Bernoulli equation
\[
\frac{dy}{dx} = \left( \eps x + b_{20} x^2 \right) \frac{1}{y} +
b_{02} y,
\]
which corresponds to cases $(b_2)$ and $(c_2)$.

\item[$(ii)$] If $\eps=b_{20}=0$ we obtain the linear equation (and, of course, integrable in a Liouville sense)
\[
\frac{dy}{dx} = \left( \lambda + b_{11}x \right) + b_{02} y.
\]
This possibility is not taken into account by Llibre and
Valls~\cite{LlibreValls} since this equation is not,
strictly speaking, in Riccati form.

\item[$(iii)$] $b_{20}=b_{02}=\lambda=0$ and $\eps b_{11} \neq 0$
(case $(a_2)$) gives rise to $dy/dx = b_{11}x + \eps x
y^{-1}$, which is a separable equation (and a Bernoulli as well) and whose
solutions are all Liouvillian.

\item[$(iv)$] If $b_{02}=0$ we obtain a Li\'enard equation,
\[
y\frac{dy}{dx} = \left( \lambda + b_{11}x
\right)y + \left( \eps x + b_{20} x^2 \right),
\]
that will be considered more deeply in a forthcoming section.

\end{itemize}

 \subsection{Families of orthogonal polynomials}

Recall that the Hypergeometric equation, including confluences, is a
particular case of the differential equation
\begin{equation}\label{eqorpol}
y'' +\frac{L}{Q}y' +\frac{\lambda}{Q} y=0,\quad
\lambda\in\mathbb{C},\quad L=a_0+a_1x,\quad Q=b_0+b_1x+b_2x^2.
\end{equation}
It is well known (see, for example,~\cite{ch}) that classical orthogonal and
Bessel polynomials are solutions of equation~\eqref{eqorpol} for suitable values of $a_j$, $b_j$ and $\lambda$. Namely,
\begin{itemize}
\item Hermite $H_n$,
\item Chebyshev of first kind $T_n$,
\item Chebyshev of second kind $U_n$,
\item Legendre $P_n$,
\item Laguerre $L_n$,
\item associated Laguerre $L_n^{(m)}$,
\item Gegenbauer $C_n^{(m)}$,
\item Jacobi $\mathcal{P}_n^{(m,\nu)}$
\item Bessel $B_n$,
\end{itemize}
where
\begin{center}
%\label{table}
\begin{tabular}[t]{|c|c|c|c|c|} \hline
\text{Family}& $Q$ & $L$   & $\lambda$ \\ \hline
$H_n$&$1$&$-2x$&$2n$\\ \hline
$T_n$&$1-x^2$&$-x$&$n^2$\\ \hline
$U_n$&$1-x^2$&$-3x$&$n(n+2)$\\ \hline
$P_n$&$1-x^2$&$-2x$&$n(n+1)$\\ \hline
$L_n$&$x$&$1-x$&$n$\\ \hline
$L_n^{(m)}$&$x$&$m+1-x$&$n$\\ \hline
$C_n^{(m)}$&$1-x^2$&$-(2m+1)x$&$n(n+2m)$\\ \hline
$\mathcal P_n^{(m,\nu)}$&$1-x^2$&$\nu-m-(m+\nu+2)x$&$n(n+1+m+\nu)$\\ \hline
$B_n$&$x^2$&$2(x+1)$&$-n(n+1)$\\ \hline
\end{tabular}
\end{center}
\medskip

Integrability conditions and solutions of
equation~\eqref{eqorpol} can be obtained
applying Kovacic algorithm (Case 1 of the algorithm).
Besides, they can also be achieved via Kimura and
Martinet-Ramis Theorems {and the \textit{parabolic cylinder equation}}
(see \cite{dulo,ki,marram}).

As a consequence, from this table,  we obtain the following result.
\begin{proposition} \label{orthogonalth}
We consider $Q$, $L$ and $\lambda$ as in the previous table.
Then, for any $\mu\neq 0$, the planar quadratic polynomial vector field
\begin{equation} \label{lien}
 %\left\{
\begin{array}{rcl}
{\ds \frac{dv}{dt}} &=& \dfrac{\lambda}{\mu}Q + (Q'-L)v + \mu v^2,\\[1.2ex]
{\ds \frac{dx}{dt}} &=& Q,
\end{array}
%\right.
\end{equation}
has invariant algebraic curves of the form ${\ds \mu v + \frac{Q(x)\mathcal{P}'_n(x)}{\mathcal{P}_n(x)}=0 }$ where $\mathcal{P}_n(x)$ is any orthogonal polynomial associated to $\lambda,Q,L$ and $n\in\N.$
\end{proposition}
\proof
The associated (Riccati) foliation of system \eqref{lien} is written into the form
\begin{equation}\label{folvect}
\frac{d {v}}{d x} = \frac{\lambda}{\mu} + \frac{Q'-L}{Q} {v} + \frac{\mu}{Q} {v}^2.
\end{equation}
Performing the change $\tilde{v}=\mu v$ equation \eqref{folvect} becomes
\[
\frac{d \tilde{v}}{d x} = \lambda + \frac{Q'-L}{Q} \tilde{v} + \frac{1}{Q} \tilde{v}^2,
\]
and can  be led into the form~\eqref{eqorpol} through the transformation $\tilde{v}=-Qy'/y$ (see transformation $\Bb$ of
Proposition~\ref{prop:almp}). For fixed $\lambda, Q$ and $L$, let $\mathcal{P}_n$ be any orthogonal polynomial solution of equation~\eqref{eqorpol}. 
Then $v_n=-(Q\mathcal{P}'_n) / (\mu\mathcal{P}_n)$ is a rational solution of  equation \eqref{folvect}. Therefore the curve $v-v_n=0$ is an invariant curve 
of the vector field \eqref{lien} for any  $n\in\N$.
\qed

According to Proposition~\ref{prop:almp} equation \eqref{folvect} can be reduced to the form $\xi''=\rho \xi$ with
\[
\rho= \frac{1}{2} \left( \frac{L}{Q} \right)' - \frac{\lambda}{Q} + \left( \frac{L}{2Q} \right)^2,
\]
and  $\xi = \mathcal{P}_n e^{ \int \frac{L}{2Q}} $ is a solution for any $n\in \N$.
 Additionally, we notice that we fall in Case~1 of Kovacic algorithm.

\subsection{Li\'enard equation}

Let us consider first order differential equations whose associated foliation
can be expressed into the Li\'enard form
\begin{equation}
y y' = f(x)y+g(x), \label{LE}
\end{equation}
with $y=y(x)$ and rational functions $f(x)$ and $g(x)$. We are concerned with the problem of
obtaining criteria on $f(x)$ and $g(x)$ such that equation~\eqref{LE} can be led into a
Riccati equation.

This is a difficult problem and, as far as the authors know, only partial answers have been given to it.
In what follows we give some examples of such results coming from the handbook~\cite{PZ} and papers~\cite{CHT1,CHT2}.

A first example is given by the $5$-parametric family (1.3.3.11 in \cite{PZ})
\begin{equation}
y y' = (a(2m+k)x^{2k}+b(2m-k)x^{m-k-1})y-(a^2mx^{4k}+cx^{2k}+b^2m)x^{2m-2k-1}, \label{LE1}
\end{equation}
$a,b,c,m,k$ being complex parameters. In order that (\ref{LE1})  come from a polynomial vector field we have that
$m$ and $k$ must be  rational numbers (see~\cite{PZ}).
The change $w=x^k$, $y=x^m(z+ax^k+bx^{-k})$ leads~\eqref{LE1} into the Riccati form
\begin{equation}
(-mz^2+2abm-c)w^\prime(z) = bk+kzw+akw^2, \label{LE2}
\end{equation}
whose associated second order linear equation is a Riemann equation. More precisely, by Lemma~\ref{prop:almp},
%({\bf J: aqu� hay que escribir expl�citamente el cambio})
it can be written as a Legendre equation
\begin{equation}
(1-t^2)u^{\prime\prime}(t) - 2tu^\prime (t)+\left(\nu (\nu+1)-\frac {\mu^2}{1-t^2}\right)u(t)=0, \label{LE3}
\end{equation}
with
\[
\mu=-\frac {m+k}{2m},
\]
and $\nu$ being a solution of
\[
\nu^2+\nu+\frac{m^2-k^2}{4m^2} -\frac {abk^2}{mc-2abm^2}=0.
\]
The difference of exponents in (\ref{LE3}) is $\mu$, $\mu$ and
$2\nu-1$ and, therefore, we are under the hypotheses of Kimura's Theorem
(see Appendix~\ref{app:se:Kimura}).
\begin{proposition}
Legendre equation~\eqref{LE3} is integrable if and only if,
either
\begin{enumerate}
\item $\mu\pm \nu\in\mathbb{Z}$ or $\nu\in\mathbb{Z}$, or
\item $\pm \mu$, $\pm \mu$, $\pm (2\nu+1)$ belong to one of the following seven families
\[
\begin{array}{|c|c|c|c|} \hline
\mbox{Case} & \mu \in & \nu \in & \mu + \nu \in \\ \hline
(a) & \Z +\frac{1}{2}  & \C                          &                 \\[1.2ex]  \hline
(b) & \Z\pm\frac{1}{3} & \frac{1}{2}\Z\pm\frac{1}{3} & \Z+\frac{1}{6} \\[1.2ex]  \hline
(c) & \Z\pm\frac{2}{5} & \frac{1}{2}\Z\pm\frac{1}{5} & \Z+\frac{1}{0} \\[1.2ex]  \hline
(d) & \Z\pm\frac{1}{3} & \frac{1}{2}\Z\pm\frac{2}{5} & \Z+\frac{1}{10} \\[1.2ex]  \hline
(e) & \Z\pm\frac{1}{5} & \frac{1}{2}\Z\pm\frac{2}{5} & \Z+\frac{1}{10} \\[1.2ex] \hline
(f) & \Z\pm\frac{2}{5} & \frac{1}{2}\Z\pm\frac{1}{3} & \Z+\frac{1}{6} \\[1.2ex] \hline
\end{array}
\]

%\begin{enumerate}
%\item $\mu\in\mathbb{Z}+\frac12$, $\nu\in\mathbb{C}$
%\item $\mu\in\mathbb{Z}\pm\frac13$, $\nu\in \frac12\mathbb{Z}\pm\frac13$
%and $\mu+ \nu\in\mathbb{Z}+\frac16$
%\item $\mu\in\mathbb{Z}\pm\frac25$,
%$\nu\in\frac12\mathbb{Z}\pm\frac{1}{5}$ and $\mu+ \nu\in
%\mathbb{Z}+\frac1{10}$
%\item $\mu\in\mathbb{Z}\pm\frac13$,
%$\nu\in\frac12\mathbb{Z}\pm\frac{2}{5}$ and $\mu+ \nu\in
%\mathbb{Z}+\frac1{10}$
%\item $\mu\in\mathbb{Z}\pm\frac15$,
%$\nu\in\frac12\mathbb{Z}\pm\frac{2}{5}$ and $\mu+ \nu\in
%\mathbb{Z}+\frac1{10}$
%\item $\mu\in\mathbb{Z}\pm\frac25$,
%$\nu\in\frac12\mathbb{Z}\pm\frac{1}{3}$ and $\mu+ \nu\in
%\mathbb{Z}+\frac16$
%\end{enumerate}
\end{enumerate}
\end{proposition}

\proof In Kimura's Theorem, the difference of exponents $\mu$, $\mu$ and
$2\nu+1$ correspond to the possibilities listed above. Indeed, they
are cases $(i)$, $(ii.1)$, $(ii.3)$, $(ii.11)$, $(ii.12)$, $(ii.13)$ and $(ii.15)$ of Kimura's
table (see Appendix~\ref{app:se:Kimura}). For the case $(i)$ we have
\begin{itemize}
\item $\mu+\mu+2\nu+1\in 2\mathbb{Z}+1\Rightarrow \mu+\nu\in\mathbb{Z}$,
\item $-\mu+\mu+2\nu+1\in 2\mathbb{Z}+1\Rightarrow \nu\in\mathbb{Z}$,
\item $\mu-\mu+2\nu+1\in 2\mathbb{Z}+1\Rightarrow \nu\in\mathbb{Z}$,
\item $\mu+\mu-2\nu-1\in 2\mathbb{Z}+1\Rightarrow \mu-\nu\in\mathbb{Z}$.
\end{itemize}
The rest of the cases can be proven in a similar way.
\begin{itemize}
\item[(ii.1)] We see that $\pm \mu\in\frac12+\mathbb{Z}$ and $\pm(2\nu+1)\in \mathbb{C}$ and therefore $\mu\in\mathbb{Z}+\frac12$ and $\nu\in\mathbb{C}$.
\item[(ii.3)] We consider that $\pm \mu=l+\frac23$, $\pm \mu=m+\frac13$, $\pm(2\nu+1)=q+\frac13$, being $l,m,q\in\mathbb{Z}$. Take for instance $\mu\in\mathbb{Z}\pm\frac13$ and $\nu\in\frac12\mathbb{Z}\pm\frac13$. Furthermore, $l+m+q$ must be even and in consequence we obtain $\mu+\nu\in \mathbb{Z}+\frac16$.
\item[(ii.11)]We have that $\pm \mu=l+\frac25$, $\pm \mu=m+\frac25$, $\pm(2\nu+1)=q+\frac25$ with $l,m,q\in\mathbb{Z}$. For example we consider $\mu\in\mathbb{Z}\pm\frac25$ and $\nu\in\frac12\mathbb{Z}\pm\frac15$. Moreover $l+m+q$ must be even and therefore we have that  $\mu+\nu\in \mathbb{Z}+\frac{1}{10}$.
\item[(ii.12)]Now we consider $\pm \mu=l+\frac23$, $\pm \mu=m+\frac13$, $\pm(2\nu+1)=q+\frac15$, being $l,m,q\in\mathbb{Z}$. We take for instance $\mu\in\mathbb{Z}\pm\frac13$ and $\nu\in\frac12\mathbb{Z}\pm\frac25$. Additionally $l+m+q$ must be even and so  $\mu+\nu\in \mathbb{Z}+\frac{1}{10}$.
\item[(ii.13)] Let be $\pm \mu=l+\frac23$, $\pm \mu=m+\frac13$, $\pm(2\nu+1)=q+\frac15$ with $l,m,q\in\mathbb{Z}$. Consider  for example $\mu\in\mathbb{Z}\pm\frac15$ and $\nu\in\frac12\mathbb{Z}\pm\frac25$. Furthermore, $l+m+q$ must be even and therefore  $\mu+\nu\in \mathbb{Z}+\frac{1}{10}$.
\item[(ii.15)] Consider $\pm \mu=l+\frac35$, $\pm \mu=m+\frac25$, $\pm(2\nu+1)=q+\frac13$, being $l,m,q\in\mathbb{Z}$, take for instance $\mu\in\mathbb{Z}\pm\frac25$ and $\nu\in\frac12\mathbb{Z}\pm\frac13$. Moreover, $l+m+q$ must be even and in consequence  $\mu+\nu\in \mathbb{Z}+\frac{1}{6}$.
\end{itemize}
Finally, observe that the difference exponents $\mu$, $\mu$ and $2\nu+1$ do not satisfy the conditions $(ii.2)$, $(ii.4)$, $(ii.5)$, $(ii.6)$, $(ii.7)$, $(ii.8)$, $(ii.9)$, $(ii.10)$
and $(ii.14)$ %by contradiction over $\mu$. 
\qed

Now we deal with an example from \cite{CHT1,CHT2}. We consider the equation
\begin{equation}
 \frac{dx}{dw}=A(x)+B(x)w,
 \label{eq1}
 \end{equation}
where $x$ is the dependent variable and $w$ is the independent one.
Now, for $B\not\equiv 0$, by means of the change of variable
\[
w=y-\frac{A}{B},
\]
and changing $(w,x)\longrightarrow (x,y)$ (that is, we consider now $x$ as the independent variable and $y$ as the dependent one)
equation \eqref{eq1} is transformed into the Li\'enard equation
\begin{equation}
y\frac{dy}{dx}=\frac{1}{B}+\frac{d}{dx}\left( \frac{A}{B} \right)y,
\label{l2}
\end{equation}
for any functions $A$ and $B \not\equiv 0$.
In particular, for
\[
A=A(x)=a+bx+cx^2,\quad B=B(x)=\alpha+\beta x+\gamma x^2,
\]
equation \eqref{eq1} falls into the Riccatti form
%\eqref{lemma:riccati:initial_form}
\[
\dfrac{dx}{dw}=(a+\alpha w)+(b+\beta w)x+(c+\gamma w)x^2.
\]

 %where $a_0=c+\gamma x$, $a_1=b+\beta x$ and $a_2=a+\alpha x$.

By Proposition~\ref{trns}, applying the transformation $\mathcal{T}$ it follows the reduced Riccatti equation \eqref{lemma:riccati:normal_form} and through the transformation $\mathcal{R}$
the normalized second order differential equation $\xi''=\rho(x) \xi$
with
\begin{equation}
\rho(x)= \frac{\beta^2 - 4\alpha\gamma}{4} x^2 - \frac{2a\gamma + 2\alpha c - b\beta}{2} x - \frac{4ac - b^2}{4} +
\frac{b\gamma - \beta c}{2(\gamma x + c)} + \frac{3\gamma^2}{4(\gamma x + c)^2}.
\end{equation}
Introducing the change $\tau=\gamma x + c$ we get
\begin{equation}\label{eqrrcht2}
\xi''=\rho(\tau)\xi,
\end{equation}
where
\begin{eqnarray*}
\rho(\tau) &=& \frac{\beta^2 - 4\alpha\gamma}{4\gamma^2} \tau^2 - \frac{2a\gamma^2 - 2\alpha c\gamma - \beta b\gamma + \beta^2 c}{2\gamma^2}\tau +\frac{b^2\gamma^2 - 2b\beta c\gamma + \beta^2c^2}{4\gamma^2} +  \\
&& \frac{b\gamma - \beta c}{2\tau} + \frac{3\gamma^2}{4\tau^2},
\end{eqnarray*}
and performing $z=\sqrt[4]{\frac{\beta^2 - 4\alpha\gamma}{4\gamma^2} } \tau$ we arrive to
\begin{equation}
\label{eqrrcht3}
\psi''=\phi(z)\psi, \quad \phi(z)=z^2+\delta_1z+ \frac{\delta_1^2}{4} -\delta_2+ \frac{\delta_3}{2z} + \frac{\delta_0^2-1}{4z^2},
\end{equation}
with $\delta_i$ being algebraic functions in $a$, $b$, $c$, $\alpha$, $\beta$ and $\gamma$. This equation \eqref{eqrrcht3} is exactly the \emph{biconfluent Heun equation} whose integrability is analysed in Appendix B.2.

Assuming $\beta=\gamma=0$ we obtain a Li\'enard equation which is
transformable into a reduced second order differential equation with
$r\in\mathbb{C}[x]$ and $\mathrm{deg}(r)=1$. This means that the
equation is not integrable (see \cite{ko} and
section~\ref{se:other_families}). As a particular case, we have a Li\'enard equation that
can be reduced to the Riccati equation given in \cite[equation 1.3.2.1]{PL},
\begin{equation}
2yy^\prime= (ax+b)y+1.
\label{LE5}
\end{equation}

Moreover, under some
restrictions over the parameters of the biconfluent Heun equation
one can obtain the Whittaker equation. For instance, the Li\'enard
equation
\[
y\frac{dy}{dx} = \left( \lambda + b_{11}x \right)y + \left( \eps x +
b_{20} x^2 \right),
\]
falls into a Whittaker equation for some special values of the
parameters.

\begin{remark} We would like to stress the following facts.
\begin{itemize}
\item[(a)]
It is well--known that via the change $z(x)=\int f(x) dx$ (with
inverse $x=x(z)$), the Li\'enard equation (\ref{LE}) can be
led into the equation
\begin{equation}
y y' (z) = y+h(z),
\end{equation}
with
\[
h(z):=\frac {g(x(z))}{f(x(z))}.
\]
In a similar way, the change $z(x)=\int g(x) dx$ reduces~\eqref{LE} to
\begin{equation}
y y' (z) = h(z)y+1,
\end{equation}
with
\[
h(z):=\frac {f(x(z))}{g(x(z))}.
\]
However, in general, these transformations do not preserve the differential field of the coefficients. This is why we do not use them to reduce equation~\eqref{LE}.

\item[(b)] Sometimes equation (\ref{LE}) is called Abel equation of second
kind since through the change $y=1/w$ it is reduced to an
Abel equation
\[
w'=-f(x)w^2-g(x)w^3.
\]
\end{itemize}
\end{remark}
In Cheb--Terrab \cite{CHT1,CHT2}  the authors consider some families of Abel equations which are reducible to Riccati equations.

\subsection{Other families}
\label{se:other_families}
Here we consider some special cases of Riccati equations.
\medskip

\noindent a) {\it Polynomial Riccati equations}\\
The Riccati equation
\[
w'=r(x)-w^2,\quad r(x)\in\mathbb{C}[x]
\]
has been studied by several authors (see, for instance~\cite{acbl,Zol,Vi}). The Galois group of its associated second order linear differential equation is connected and can be either $\mathrm{SL}(2,\mathbb{C})$ or the Borel group~\cite{acmowe,acbl}. In the first case the tangent field associated to its Riccati equation has no invariant curves. In the second one, there is no rational first integral for its vector field.
As an example, the reduced form for the \emph{triconfluent Heun equation} is of this type and is given by
\begin{equation}\label{tricheun}
\xi''=\rho(x)\xi,\quad \rho(x)=\frac{9x^4}{4}+\frac{3}{2}\delta_2 x^2-\delta_1 x+\frac{\delta_2^2}{4}-\delta_0.
\end{equation}

\noindent b) {\it Equations with finite Galois group}\\
Consider the  polynomial Riccati vector field
\begin{equation}\label{odeweil0}
\begin{array}{lll}\dot x&=&-x^2(x-1)^2\\ \dot y&=& a(x-1)^2+bx^2+cx(x-1)+x^2(x-1)^2v^2.
\end{array}
\end{equation}
Its foliation was studied in~\cite{weil2} and is related to the following differential equation
\begin{equation}\label{odeweil1}
\xi''=- \left(\frac {1-\lambda^2}{4x^2}+\frac{1-\mu^2}{4(x-1)^2}+\frac{1-\nu^2+\lambda^2+\mu^2}{4x(x-1
) } \right)\xi,
\end{equation}
where $\lambda$, $\mu$ and $\nu$ are the differences of the exponents at $0$, $1$ and $\infty$ of the hypergeometric equation, see also Appendix \ref{appendix:B}.
 When the equation \eqref{odeweil1} is integrable their solutions are Legendre functions. For  $(\lambda,\mu,\nu)= (\frac12, \frac12, \frac{1}{n})$,
the solutions of  equation~\eqref{odeweil1} are given by
\begin{eqnarray*}
\xi_1 &=& \sqrt [4]{x^2-x} \left( 2x-1+2\sqrt {x^2-x}
 \right) ^{\frac12\sqrt {\frac {1-2n^2}{n^2}}},\\
\xi_2 &=& \sqrt [4]{x^2-x} \left( 2x-1+2\sqrt {x^2-x}
 \right) ^{-\frac12\sqrt {\frac {1-2n^2}{n^2}}},
\end{eqnarray*}
and its differential Galois group is the \emph{Dihedral Group} $D_n$ (case 2 of Kovacic algorithm). One can also obtain the tetrahedral,
 octahedral and icosahedral groups for $(\lambda,\mu,\nu)=(\frac13,\frac12,\frac13)$, $(\lambda,\mu,\nu)=(\frac13,\frac12,\frac14)$
and $(\lambda,\mu,\nu)=(\frac13,\frac12,\frac15)$, respectively (case 3 of Kovacic algorithm).

Related to~\eqref{odeweil1} we have the equation (see \cite{weil2})
\begin{equation}\label{odeweil2}
y''+\frac{7x-4}{6x( x-1)}y'
-\frac{36\nu^2-1}{144x(x-1)}y=0,
\end{equation}
where the differential Galois group can be tetrahedral, octahedral and icosahedral depending on the values of
 $\nu$ (indeed, for $1/3,1/4$ and $1/5$ respectively). Even though explicit solutions of equations \eqref{odeweil1} and \eqref{odeweil2}
are difficult to get in general, it can be proved the existence of rational first integrals for the associated vector field for values of the
parameters $\nu\in\{\frac13,\frac14,\frac15\}$. For instance, solutions of equation \eqref{odeweil2} for $\nu=\frac13$ (tetrahedral group) are
\[
\begin{array}{lll}
y_1 &=&\sqrt [4]{-\sqrt [3]{x}+1}\sqrt [8]{{\frac {2\sqrt {\sqrt[3]{x^2}+\sqrt [3]{x}+1}-\sqrt {3}\sqrt [3]{x}-\sqrt {3}}{
-\sqrt {3}\sqrt [3]{x}-\sqrt {3}-2\sqrt {\sqrt[3]{x^2}+
\sqrt [3]{x}+1}}}},\\\\
y_2 &=& {{\sqrt [4]{-\sqrt [3]{x}+1}}{\sqrt [8]{{\frac {2\sqrt {\sqrt[3]{x^2}+\sqrt [3]{x}+1}-\sqrt {3}\sqrt [3]{x}-\sqrt {3}}{
-\sqrt {3}\sqrt [3]{x}-\sqrt {3}-2\sqrt {\sqrt[3]{x^2}+
\sqrt [3]{x}+1}}}}}}.
\end{array}
\]
In the case $\nu=\frac14$ (octahedral group) we obtain the solutions
${\ds y_1=\frac{f_1 f_2}{f_3}}$ and ${\ds y_2=\frac{f_1}{f_2f_3} }$ where
\[
\begin{array}{lll}
f_1&=&\sqrt [6]{\sqrt [3]{{\frac {\sqrt {1-x}+1}{\sqrt {1-x}-1}}}-1},\\
f_2&=&\sqrt [8]{-\frac{-2\,\sqrt { \left( {\frac {\sqrt {1-x}+1}{\sqrt {1-x}-1}} \right) ^{2/
3}+\sqrt [3]{{\frac {\sqrt {1-x}+1}{\sqrt {1-x}-1}}}+1}+\sqrt {3}
\sqrt [3]{{\frac {\sqrt {1-x}+1}{\sqrt {1-x}-1}}}+\sqrt {3}
}{\sqrt {3}\sqrt [3]{{\frac {\sqrt {1-x}+1}{\sqrt {1-x}-1}}}+\sqrt {3}+2
\,\sqrt { \left( {\frac {\sqrt {1-x}+1}{\sqrt {1-x}-1}} \right) ^{2/3}
+\sqrt [3]{{\frac {\sqrt {1-x}+1}{\sqrt {1-x}-1}}}+1}
}},\\
f_3&=& \left(  \left( {\frac {\sqrt {1-x}+1}{\sqrt {1-x}-1}} \right) ^{2/3}+
\sqrt [3]{{\frac {\sqrt {1-x}+1}{\sqrt {1-x}-1}}}+1 \right) ^{\frac{1}{12}}.
\end{array}
\]
Finally, setting $\nu=\frac15$ (icosahedral group), we get the solutions
${\ds y_1=e^{\int\omega_1}}$ and ${\ds y_2=e^{\int\omega_2} }$,
where $\omega_1=\omega_1(x)$ and $\omega_2=\omega_2(x)$ are roots of the polynomial ${\ds \sum_{k=0}^{12}a_k(x)\omega^{k}}$ and
\begin{eqnarray*}
a_0(x) &=& 102400-11264x-11x^2\\
a_1(x) &=& -3686400x+3679200x^2+7200x^3\\
a_2(x) &=& 479001600x^2-476863200x^3-2138400x^4\\
a_3(x) &=& 30412800000x^2 -60445440000x^3 + 29652480000x^4 + 380160000x^5\\
a_4(x) &=& -821145600000x^3+1597384800000x^4 -731332800000x^5 -44906400000x^6\\
a_5(x) &=& -3695155200000x^4+11085465600000x^5+3695155200000x^7-\\
&& 11085465600000x^6\\
\end{eqnarray*}
\begin{eqnarray*}
a_6(x) &=& -492687360000000x^4 + +1693612800000000x^5 -2124714240000000x^6 + \\
&& 1139339520000000x^7-215550720000000x^8\\
a_7(x) &=& 8868372480000000x^5-35473489920000000x^6 + 53210234880000000x^7 - \\
&& 35473489920000000x^8 + 8868372480000000x^9 ,\\
a_8(x) &=& -249422976000000000x^6+997691904000000000x^7-1496537856000000000x^8 +\\
&& 997691904000000000x^9 -249422976000000000x^{10}\\
a_9(x) &=&  -4434186240000000000x^6 + 22170931200000000000x^7 - \\
&& 44341862400000000000x^8 + 44341862400000000000x^9 - \\
&& 22170931200000000000x^{10} + 4434186240000000000x^{11}\\
a_{10}(x) &=& 39907676160000000000x^7-199538380800000000000x^8+ \\
&& 399076761600000000000x^9 -399076761600000000000x^{10} + \\
&& 199538380800000000000x^{11} -39907676160000000000x^{12}\\
a_{11}(x) &=& 0\\
a_{12}(x) &=& 2176782336000000000000x^8-13060694016000000000000x^9 + \\
&& 32651735040000000000000x^{10} -43535646720000000000000x^{11}+ \\
&& 32651735040000000000000x^{12} -13060694016000000000000x^{13} + \\ && 2176782336000000000000x^{14}.
\end{eqnarray*}

\noindent c) {\it Lam\'e families}.\\

Let us consider in the Lam\'e equation
\begin{equation}
\frac{d^2 y}{dx^2}+\frac{f^{\prime}(x)}{2f(x)}\frac{d y}{dx}-
\frac {n(n+1)x + B} {f(x)}  y = 0, \label{LAM}
\end{equation}
the case where $f(x)=4x^3-g_2x-g_3$ and parameters $n$, $B$, $g_2$, $g_3$
such that the discriminant of $f$, namely $27g_3^2-g_2^3$, is
non-zero (see Appendix~\ref{appendix:B}). Performing the transformation
\begin{equation}v=-\frac {y^\prime}{cy},
 \label{CR}
\end{equation}
with $c=c(x)$ being any non-zero arbitrary rational function, it follows the family of Riccati equations associated to the Lam\'e equation,
\begin{equation}v^\prime= -
\frac {n(n+1)x + B} {cf}-\left(\frac{f^{\prime}(x)}{2f(x)}+\frac{c^\prime}{c}\right)v+cv^2.
 \label{RILA}
\end{equation}

For the Lam\'e case, (i.1) of \ref{seclame} with $B=B_i$, and from
Remark \ref{remlamesol} it turns out the existence of a polynomial solution of equation
\eqref{RILA} despite of the general solution of this equation is not
algebraic. Furthermore, it is clear that for a fixed $n$ we have
associated Riccati equations (\ref{RILA}) with arbitrary degree,
because the non-zero rational function $c=c(x)$ at the transformation above is arbitrary. The Lam\'e functions correspond here to algebraic solutions of (\ref{RILA}) and  moving $n$ along the integers, we  obtain
algebraic invariants curves
\[
v+\frac {E^\prime(x)}{c(x)E(x)}=0,
\]
of unbounded degree.

Hence, we have prove the following result.
\begin{theorem}\label{teoremita1}
For any fixed degree in the Riccati family  (\ref{RILA}) associated to the Lam\'e equation (\ref{LAM}) there exist invariant algebraic curves
of any tangent field $X$ of the corresponding foliation with unbounded degree. Furthermore, the first integral of $X$ is not rational.
\end{theorem}

%SSSSSSSSSSSSSSSSSSSSSSSSSSSSSSSSSSSSSSSSSSSSSSSSSSSSSSSSSSSSSSSSSSSSSSSSSSSSSSs

\section*{Appendix}\label{appendix}

\appendix
%\section{Algorithmic Considerations}

\section{Kovacic Algorithm}\label{kovsection}
This algorithm is devoted to solve the reduced linear differential equation (RLDE) $\xi''=\rho \xi$ and is
based on the algebraic subgroups of $\mathrm{SL}(2,\mathbb{C}).$ For more
details see \cite{ko}. Although improvements for this algorithm are given
in \cite{dulo,ulwe}, we follow the original version given by Kovacic in
\cite{ko}.

\begin{theorem}\label{subgroups} Let $G$ be an algebraic subgroup of $\mathrm{SL}(2,\mathbb{C})$.  Then
one of the following four cases can occur.
\begin{enumerate}
\item $G$ is triangularizable.
\item $G$ is conjugate to a subgroup of infinite dihedral group (also called meta-abelian group)
and case 1 does not hold.
\item Up to conjugation $G$ is one of the following finite groups: Tetrahedral group, Octahedral
group or Icosahedral group, and cases 1 and 2 do not hold.
\item $G = \mathrm{SL}(2,\mathbb{C})$.
\end{enumerate}
\end{theorem}

Each case in Kovacic algorithm is related to each one of the
algebraic subgroups of $\mathrm{SL}(2,\mathbb{C})$ and its associated
Riccatti equation
\[
\theta^{\prime}=r-\theta ^{2}=\left( \sqrt{r}-\theta\right)
\left(  \sqrt{r}+\theta\right),\quad\theta=\frac{\xi'}{\xi}, \ r=\rho.
\]
According to Theorem \ref{subgroups} we obtain four cases. Only for cases 1, 2 and 3 one can solve the
differential equation RLDE and in case 4 one has not
Liouvillian solutions for it. Kovacic
algorithm can possibly provide one solution ($\xi_1$), so
the second one ($\xi_2$) can be got through
\begin{equation}\label{second}
\xi_2=\xi_1\int\frac{dx}{\xi_1^2}.
\end{equation}

\noindent \textbf{Notations.} For the RLDE given by
\[
\xi''=\rho\xi=r\xi,\qquad r=\frac{s}{t},\quad s,t\in \mathbb{C}[x],
\]
we use:
\begin{enumerate}
\item Denote by $\Gamma'$ be
the
set of (finite) poles of $r$, $\Gamma^{\prime}=\left\{  c\in\mathbb{C}%
:t(c)=0\right\}$.

\item Denote by
$\Gamma=\Gamma^{\prime}\cup\{\infty\}$.
\item By the order of $r$ at
$c\in \Gamma'$, $\circ(r_c)$, we mean the multiplicity of $c$ as a
pole of $r$.

\item By the order of $r$ at $\infty$, $\circ\left(
r_{\infty}\right) ,$ we mean the order of $\infty$ as a zero of
$r$. That is $\circ\left( r_{\infty }\right) =deg(t)-deg(s)$.

\end{enumerate}
\subsection{The four cases}

\noindent\textbf{Case 1.} In this case $\left[ \sqrt{r}\right] _{c}$ and
$\left[ \sqrt{r}\right] _{\infty}$ means the Laurent series of
$\sqrt{r}$ at $c$ and the Laurent series of $\sqrt{r}$ at $\infty$
respectively. Furthermore, we define $\varepsilon(p)$ as follows:
if $p\in\Gamma,$ then $\varepsilon\left( p\right) \in\{+,-\}.$
Finally, the complex numbers $\alpha_{c}^{+},\alpha_{c}^{-},\alpha_{\infty}%
^{+},\alpha_{\infty}^{-}$ will be defined in the first step. If
the differential equation has not poles it only can fall in this
case.
\medskip

\textbf{Step 1.} Search for each $c \in \Gamma'$ and for $\infty$ the
corresponding situation as follows:

\medskip

\begin{itemize}

\item[$(c_{0})$] If $\circ\left(  r_{c}\right)  =0$, then
$$\left[ \sqrt {r}\right] _{c}=0,\quad\alpha_{c}^{\pm}=0.$$

\item[$(c_{1})$] If $\circ\left(  r_{c}\right)  =1$, then
$$\left[ \sqrt {r}\right] _{c}=0,\quad\alpha_{c}^{\pm}=1.$$

\item[$(c_{2})$] If $\circ\left(  r_{c}\right)  =2,$ and $$r= \cdots
+ b(x-c)^{-2}+\cdots,\quad \text{then}$$
$$\left[ \sqrt {r}\right]_{c}=0,\quad \alpha_{c}^{\pm}=\frac{1\pm\sqrt{1+4b}}{2}.$$

\item[$(c_{3})$] If $\circ\left(  r_{c}\right)  =2v\geq4$, and $$r=
(a\left( x-c\right)  ^{-v}+...+d\left( x-c\right)
^{-2})^{2}+b(x-c)^{-(v+1)}+\cdots,\quad \text{then}$$ $$\left[
\sqrt {r}\right] _{c}=a\left( x-c\right) ^{-v}+...+d\left(
x-c\right) ^{-2},\quad\alpha_{c}^{\pm}=\frac{1}{2}\left(
\pm\frac{b}{a}+v\right).$$

\item[$(\infty_{1})$] If $\circ\left(  r_{\infty}\right)  >2$, then
$$\left[\sqrt{r}\right]  _{\infty}=0,\quad\alpha_{\infty}^{+}=0,\quad\alpha_{\infty}^{-}=1.$$

\item[$(\infty_{2})$] If $\circ\left(  r_{\infty}\right)  =2,$ and
$r= \cdots + bx^{2}+\cdots$, then $$\left[
\sqrt{r}\right]  _{\infty}=0,\quad\alpha_{\infty}^{\pm}=\frac{1\pm\sqrt{1+4b}%
}{2}.$$

\item[$(\infty_{3})$] If $\circ\left(  r_{\infty}\right) =-2v\leq0$,
and
$$r=\left( ax^{v}+...+d\right)  ^{2}+ bx^{v-1}+\cdots,\quad \text{then}$$
$$\left[  \sqrt{r}\right]  _{\infty}=ax^{v}+...+d,\quad
\text{and}\quad \alpha_{\infty}^{\pm }=\frac{1}{2}\left(
\pm\frac{b}{a}-v\right).$$
\end{itemize}
\medskip

\textbf{Step 2.} Find $D\neq\emptyset$ defined by
$$D=\left\{
m\in\mathbb{Z}_{+}:m=\alpha_{\infty}^{\varepsilon
(\infty)}-%
{\displaystyle\sum\limits_{c\in\Gamma^{\prime}}}
\alpha_{c}^{\varepsilon(c)},\forall\left(  \varepsilon\left(
p\right) \right)  _{p\in\Gamma}\right\}  .$$ If $D=\emptyset$,
then we should start with the case 2. Now, if $\#D>0$, then for
each $m\in D$ we search $\omega$ $\in\mathbb{C}(x)$ such that
$$\omega=\varepsilon\left(
\infty\right)  \left[  \sqrt{r}\right]  _{\infty}+%
{\displaystyle\sum\limits_{c\in\Gamma^{\prime}}}
\left(  \varepsilon\left(  c\right)  \left[  \sqrt{r}\right]  _{c}%
+{\alpha_{c}^{\varepsilon(c)}}{(x-c)^{-1}}\right).$$
\medskip

\textbf{Step 3}. For each $m\in D$, search for a monic polynomial
$P_m$ of degree $m$ with
$$P_m'' + 2\omega P_m' + (\omega' + \omega^2 - r) P_m = 0.$$

If one successes then $\xi_1=P_m e^{\int\omega}$ is a
solution of the differential equation  RLDE.  Else, Case 1
cannot hold.
\bigskip

\noindent\textbf{Case 2.}  Search for each $c \in \Gamma'$ and for
$\infty$ the corresponding situation as follows.
\medskip

\textbf{Step 1.} Search for each $c\in\Gamma'$ and $\infty$
the sets $E_{c}\neq\emptyset$ and $E_{\infty}\neq\emptyset.$ For
each $c\in\Gamma'$ and for $\infty$ we define
$E_{c}\subset\mathbb{Z}$ and $E_{\infty}\subset\mathbb{Z}$ as
follows:
\medskip

\begin{itemize}
\item[($c_1$)] If $\circ\left(  r_{c}\right)=1$, then $E_{c}=\{4\}.$

\item[($c_2$)] If $\circ\left(  r_{c}\right)  =2,$ and $r= \cdots +
b(x-c)^{-2}+\cdots ,\ $ then $$E_{c}=\left\{
2+k\sqrt{1+4b}:k=0,\pm2\right\}.$$

\item[($c_3$)] If $\circ\left(  r_{c}\right)  =v>2$, then $E_{c}=\{v\}.$

\item[$(\infty_{1})$] If $\circ\left(  r_{\infty}\right)  >2$, then
$E_{\infty }=\{0,2,4\}.$

\item[$(\infty_{2})$] If $\circ\left(  r_{\infty}\right)  =2,$ and
$r= \cdots + bx^{2}+\cdots$, then $$E_{\infty }=\left\{
2+k\sqrt{1+4b}:k=0,\pm2\right\}.$$

\item[$(\infty_{3})$] If $\circ\left(  r_{\infty}\right)  =v<2$,
then $E_{\infty }=\{v\}.$
\medskip
\end{itemize}

\textbf{Step 2.} Find $D\neq\emptyset$ defined by
$$D=\left\{
m\in\mathbb{Z}_{+}:\quad m=\frac{1}{2}\left(  e_{\infty}-
{\displaystyle\sum\limits_{c\in\Gamma^{\prime}}} e_{c}\right)
,\forall e_{p}\in E_{p},\text{ }p\in\Gamma\right\}.$$ If
$D=\emptyset,$ then we should start the case 3. Now, if $\#D>0,$
then for each $m\in D$ we search a rational function $\theta$
defined by
$$\theta=\frac{1}{2}
{\displaystyle\sum\limits_{c\in\Gamma^{\prime}}}
\frac{e_{c}}{x-c}.$$
\medskip

\textbf{Step 3.} For each $m\in D,$ search a monic polynomial $P_m$
of degree $m$, such that
$$P_m^{\prime\prime\prime}+3\theta
P_m^{\prime\prime}+(3\theta^{\prime}+3\theta
^{2}-4r)P_m^{\prime}+\left(  \theta^{\prime\prime}+3\theta\theta^{\prime}%
+\theta^{3}-4r\theta-2r^{\prime}\right)P_m=0.$$ If $P_m$ does not
exist, then Case 2 cannot hold. If such a polynomial is found, set
$\phi = \theta + P'/P$ and let $\omega$ be a solution of
$$\omega^2 + \phi \omega + \frac{1}{2} \left(\phi' + \phi^2 -2r\right)=
0.$$

Then $\xi_1 = e^{\int\omega}$ is a solution of the differential
equation  RLDE.
\bigskip

\noindent\textbf{Case 3.} Search for each $c \in \Gamma'$ and for
$\infty$ the corresponding situation as follows:
\medskip

\textbf{Step 1.} Search for each $c\in\Gamma^{\prime}$ and $\infty$
the sets $E_{c}\neq\emptyset$ and $E_{\infty}\neq\emptyset.$ For
each $c\in\Gamma^{\prime}$ and for $\infty$ we define
$E_{c}\subset\mathbb{Z}$ and $E_{\infty}\subset\mathbb{Z}$ as
follows:
\medskip

\begin{itemize}

\item[$(c_{1})$] If $\circ\left(  r_{c}\right)  =1$, then
$E_{c}=\{12\}.$

\item[$(c_{2})$] If $\circ\left(  r_{c}\right)  =2,$ and $r= \cdots +
b(x-c)^{-2}+\cdots$, then
\begin{displaymath}
E_{c}=\left\{ 6+k\sqrt{1+4b}:\quad
k=0,\pm1,\pm2,\pm3,\pm4,\pm5,\pm6\right\}.
\end{displaymath}

\item[$(\infty)$] If $\circ\left(  r_{\infty}\right)  =v\geq2,$ and $r=
\cdots + bx^{2}+\cdots$, then
$$E_{\infty }=\left\{
6+ \frac{12k}{n} \sqrt{1+4b}:\quad
k=0,\pm1,\pm2,\pm3,\pm4,\pm5,\pm6\right\},\quad n\in\{4,6,12\}.$$
\medskip
\end{itemize}

\textbf{Step 2.} Find $D\neq\emptyset$ defined by
$$D=\left\{
m\in\mathbb{Z}_{+}:\quad m=\frac{n}{12}\left(
e_{\infty}-{\displaystyle\sum\limits_{c\in\Gamma^{\prime}}}
e_{c}\right)  ,\forall e_{p}\in E_{p},\text{
}p\in\Gamma\right\}.$$ In this case we start with $n=4$ to obtain
the solution, afterwards $n=6$ and finally $n=12$. If
$D=\emptyset$, then the differential equation has not Liouvillian
solution because it falls in the case 4. Now, if $\#D>0,$ then for
each $m\in D$ with its respective $n$, search a rational function
$$
\theta=\frac{n}{12} {\displaystyle\sum\limits_{c\in\Gamma^{\prime}}}
\frac{e_{c}}{x-c},
$$
and a polynomial $S$ defined as $$S=
{\displaystyle\prod\limits_{c\in\Gamma^{\prime}}} (x-c).$$

\textbf{Step 3}. Search for each $m\in D$, with its respective $n$, a
monic polynomial $P_m=P$ of degree $m,$ such that its coefficients
can be determined recursively by
$$\bigskip P_{-1}=0,\quad P_{n}=-P,$$
$$P_{i-1}=-SP_{i}^{\prime}-\left( \left( n-i\right)
S^{\prime}-S\theta\right)  P_{i}-\left( n-i\right)  \left(
i+1\right)  S^{2}rP_{i+1},$$ where $i\in\{0,1\ldots,n-1,n\}.$ If
$P$ does not exist, then the differential equation has not
Liouvillian solution because it falls in Case 4. Now, if $P$
exists search $\omega$ such that $$
{\displaystyle\sum\limits_{i=0}^{n}} \frac{S^{i}P}{\left(
n-i\right)  !}\omega^{i}=0,$$ then a solution of the differential
equation the RLDE is given by $$\xi=e^{\int \omega},$$ where
$\omega$ is solution of the previous polynomial of degree $n$.
\bigskip

%SSSSSSSSSSSSSSSSSSSSSSSSSSSSSSSSSSSSSSSSSSSSSSSSSSSSSSSSSSSSSSSSSSSSSSSSSSSSSSSSSSSSSSSSSSSSSSSSSSSS

\section{Some Special Functions}
\label{appendix:B}

\subsection{Hypergeometric families}

\subsubsection{Kimura's Theorem}
\label{app:se:Kimura}

The hypergeometric (or Riemann) equation is the more general
second order linear differential equation over the Riemann sphere
with three regular singular singularities. If we place the
singularities at $x = 0, 1, \infty$ it is given by
\begin{eqnarray}
\label{hypergeometric1}
&& \frac{d^2 y}{dx^2}+ \left(\frac{1-\alpha-\alpha'}{x} + \frac{1-\gamma-\gamma'}{x-1}\right)\frac{d y}{dx}\\
&& \qquad\qquad\qquad \qquad  + \left(\frac{\alpha\alpha'}{x^2} +
\frac{\gamma\gamma'}{(x-1)^2} +
\frac{\beta\beta'-\alpha\alpha'\gamma\gamma'}{x(x-1)}\right) y =
0, \nonumber
\end{eqnarray}
where $(\alpha , \alpha ')$, $(\gamma , \gamma ')$, $(\beta ,
\beta ')$ are the exponents at the singular points and must
satisfy the Fuchs relation $\alpha + \alpha' + \gamma +
\gamma'+\beta + \beta'= 1$.\\

Now, we will briefly describe  Kimura's Theorem that provides
necessary and sufficient conditions for the integrability of the
hypergeometric equation. Let be ${\lambda} = \alpha -\alpha'$,
${\mu} = \beta-\beta'$ and ${\nu} = \gamma -\gamma'$.

\begin{theorem}[Kimura, \cite{ki}]\label{kimurath}
The hypergeometric
equation (\ref{hypergeometric1}) is integrable if and only if
either
\begin{enumerate}
\item[(i)] At least one of the four numbers ${\lambda}+{\mu}+{\nu}$,
$-{\lambda}+{\mu}+{\nu}$, ${\lambda}-{\mu}+{\nu}$,
${\lambda}+{\mu}-{\nu}$ is an odd integer, or
\item[(ii)] The numbers ${\lambda}$ or $-{\lambda}$, ${\mu}$
or $-{\mu}$ and ${\nu}$ or $-{\nu}$ belong (in an arbitrary order)
to some of the following fifteen families
\end{enumerate}
$$
\begin{array}{|c|c|c|c|c|}\hline
1 & 1/2+l & 1/2+m & \mbox{arbitrary complex number} &\\ \hline 2 &
1/2+l & 1/3+m & 1/3+q &\\ \hline 3 & 2/3+l & 1/3+m & 1/3+q &
l+m+q\mbox{ even}\\ \hline 4 & 1/2+l & 1/3+ m & 1/4+q & \\ \hline
5 & 2/3+l & 1/4+m & 1/4+q & l+m+q\mbox{ even}\\ \hline 6 & 1/2+l &
1/3+m & 1/5+q & \\ \hline 7 & 2/5+l & 1/3+m & 1/3+q & l+m+q\mbox{
even}\\ \hline 8 & 2/3+l & 1/5+m & 1/5+q & l+m+q\mbox{ even}\\
\hline 9 & 1/2+l & 2/5+m & 1/5+q & l+m+q\mbox{ even}\\ \hline 10 &
3/5+l & 1/3+m & 1/5+q & l+m+q\mbox{ even}\\ \hline 11 &2/5+l &
2/5+m & 2/5+q & l+m+q\mbox{ even}\\ \hline 12 &2/3+l & 1/3+m &
1/5+q & l+m+q\mbox{ even}\\ \hline 13 & 4/5+l & 1/5+m & 1/5+q &
l+m+q\mbox{ even}\\ \hline 14 & 1/2+l &2/5+m & 1/3+q & l+m+q\mbox{
even}\\ \hline 15 & 3/5+l & 2/5+m & 1/3+q & l+m+q\mbox{ even}\\
\hline
\end{array}
$$
Here $l,m$ and $q$ are integers.
\end{theorem}

\subsubsection{Confluent hypergeometric}

The \textit{confluent Hypergeometric equation} is a degenerate
form of the Hypergeometric differential equation where two of the
three regular singularities merge into an irregular singularity. The following are two classical
forms.
\begin{itemize}
\item \textit{Kummer's} form
\begin{equation}\label{kummer}
y''+\frac{c-x}{x} y'- \frac{a}{x} y=0, \qquad a,c\in\C
\end{equation}

\item \textit{Whittaker's} form
\begin{equation}\label{whittaker}
y''=\left(\frac{1}{4}- \frac{\kappa}{x} + \frac{4\mu^2-1}{4x^2} \right)y,
\end{equation}
\end{itemize}
where the parameters of the two equations are linked by
$\kappa=\frac{c}2-a$ and $\mu=\frac{c}{2}-\frac{1}{2}$. Furthermore,
using the expression \eqref{lemma:2ode:initial_form}, we can see that the
Whittaker's equation is the reduced form of the Kummer's equation \eqref{kummer}.
The Galoisian structure of these equations has been deeply studied
in \cite{marram,dulo}.

\begin{theorem}[Martinet \& Ramis, \cite{marram}]\label{thmarram}
The Whittaker's differential equation \eqref{whittaker} is
integrable if and only if either,
$\kappa+\mu\in\frac12+\mathbb{N}$, or
$\kappa-\mu\in\frac12+\mathbb{N}$, or
$-\kappa+\mu\in\frac12+\mathbb{N}$, or
$-\kappa-\mu\in\frac12+\mathbb{N}$.
\end{theorem}
The \textit{Bessel's equation} is a particular case of the
confluent Hypergeometric equation and is given by
\begin{equation}\label{bessel}
y''+ \frac{1}{x} y'+ \frac{x^2-n^2}{x^2} y=0.
\end{equation}
Under a suitable transformation, the
reduced form of the Bessel's equation is a particular case of the
Whittaker's equation \eqref{whittaker}.

\begin{corollary}\label{corbessel}
The Bessel's differential equation \eqref{bessel} is integrable if
and only if $n\in \frac12+\mathbb{Z}$.
\end{corollary}
\subsection{Heun's families}
The Heun's equation is the generic differential equation with four regular singular points at $0$, $1$, $c$ and $\infty$. In its reduced form, the Heun's equation is $y''=r(x)y$, where
\begin{equation}\label{genheun}
r(x)=\frac{A}{x}+\frac{B}{x-1}+\frac{C}{x-c}+\frac{D}{x^2}+\frac{E}{(x-1)^2}+\frac{F}{(x-c)^2},
\end{equation}
$$A=-\frac{\alpha\beta}{2}-\frac{\alpha\gamma}{2c}+\frac{\delta\eta h}{c},\quad B=\frac{\alpha\beta}{2}-\frac{\beta\gamma}{2(c-1)}-\frac{\delta\eta(h-1)}{c-1},$$
$$C=\frac{\alpha\gamma}{2c}+\frac{\beta\gamma}{2(c-1)}-\frac{\delta\eta(c-h)}{c(c-1)},\quad D=\frac{\alpha}{2}\left(\frac{\alpha}{2}-1\right),\quad E=\frac{\beta}{2}\left(\frac{\beta}{2}-1\right),$$
$$F=\frac{\gamma}{2}\left(\frac{\gamma}{2}-1\right),\quad\textrm{with}\quad \alpha+\beta+\gamma-\delta-\eta=1.$$

To our purposes we write the determinant $\Pi_{d+1}(a,b,u,v,\xi,w)$ as in \cite{dulo}:
{\footnotesize{$$\left|\begin{array}{ccccccc}
w&u&0&0&0&\ldots&0\\
d\xi w+1&v&2(u+b)&0&0&\ldots&0\\
0&(d-1)\xi&w+2(v+a)&3(u+2b)&0&\ldots&0\\
0&0&(d-2)\xi&w+3(v+2a)&4(u+3b)&\ldots&0\\
\vdots&&&&&\ldots&\\
0&\ldots&&\ldots&2\xi&w+(d-1)(v+(d-2)a)&d(u+(d-1)b)\\
0&\ldots&&\ldots&0&\xi&w+d(v+(d-1)a) \end{array}\right|$$}}

\subsubsection{Biconfluent Heun}

The equation
\begin{equation}
\label{biconfluenth} \xi''=\left( x^2+\delta_1x+ \frac{\delta_1^2}{4} -\delta_2+ \frac{\delta_3}{2x} + \frac{\delta_0^2-1}{4x^2} \right) \xi,
\end{equation}
is the well known \emph{biconfluent Heun equation} which has been deeply analyzed by Duval and Loday-Richaud in \cite[p. 236]{dulo}.

\begin{theorem}\cite{dulo}.\label{biconheun}
The biconfluent Heun equation \eqref{biconfluenth} has Liouvillian solutions if and only if it falls in Case 1 of Kovacic algorithm and one of the following conditions is fulfilled.
\begin{enumerate}

\item $\delta_0^2=1$, $\delta_3=0$ and $\delta_2\in 2{\Z}+1.$

\item $\delta_0^2=1$, $\delta_3\neq 0$ and $\delta_2\in 2{\Z}^*+1$ with $|\delta_2|\geq 3$, and if $\varepsilon=\textsf{sign}\, \delta_2$, then
\[
\Pi_{(|\delta_2|-1)/2}\left(0,1,2,\varepsilon\delta_1, -2\varepsilon, \varepsilon\delta_1- \frac{\delta_3}{2}\right)=0.
\]
\item $\delta_0\neq\pm1$, $\pm\delta_0\pm\delta_2\in 2{\Z}^*$ and if $\varepsilon_0,\varepsilon_\infty\in\{\pm1\}$ are such that $\varepsilon_\infty\delta_2-\varepsilon_0\delta_0=2d^*\in 2{\N}^*$ then
\[
\Pi_{d^*}\left(0,1,1+\varepsilon_0\delta_0, \varepsilon_\infty\delta_1,-2\varepsilon_\infty,\frac{1}{2}(\varepsilon_\infty
\delta_1(1+\varepsilon_0\delta_0)-\delta_3)\right)=0.
\]
\end{enumerate}

\end{theorem}

\subsection{Lam\'e equation} \label{seclame}

The algebraic form of the Lam\'e Equation is \cite{POO,WW}
\begin{equation}
\frac{d^2 y}{dx^2}+\frac{f^{\prime}(x)}{2f(x)}\frac{d y}{dx}-
\frac {n(n+1)x + B} {f(x)}  y = 0, \label{LAM1}
\end{equation}
where $f(x)=4x^3-g_2x-g_3$, with $n$, $B$, $g_2$ and  $g_3$
parameters such that the discriminant of $f$, $\Delta=27g_3^2-g_2^3$, is
non-zero and, therefore, it has no multiple roots. This equation is a Fuchsian differential equation with
four singular points over the Riemann sphere: the roots $e_1,e_2, e_3$ of $f$ and the point of the infinity.

The mutually-exclusive known cases of solutions in closed form of
the Lam\'e equation (\ref {LAM1}) are the following.
\begin{itemize}

\item[$(i)$] The \emph{Lam\'e-Hermite} case (see \cite{POO,WW}). We have
$n\in {\bf N}$ and arbitrary parameters $B,g_2$ and $g_3$.

\item[$(ii)$] The \emph{Brioschi-Halphen-Crawford} case (see~\cite{HAL,POO}). We have $n\in \N$ such that
$m:=n+\frac 12 \in\N$ and parameters $B$, $ g_2$ and $  g_3$ satisfying an algebraic condition
$Q_m( g_2/4, g_3/4,B)=0$,
where $Q_m\in \Z[ g_2/4, g_3/4,B]$ is a polynomial of degree $m$ in $B$, known as the
\emph{Brioschi determinant}.

\item[$(iii)$] The \emph{Baldassarri} case (see \cite{BAL}). One asks $n$
to satisfy that $n+\frac 12 \in\frac 13\Z\cup\frac 14\Z\cup\frac 15\Z-\Z$ besides some
additional (involved) algebraic restrictions on the other parameters.

\end{itemize}

It is possible to prove that the only integrable cases of the
Lam\'e equation are cases $(i)$--$(iii)$ above, see \cite{morales}. In
$(ii)$ and $(iii)$ the general solution of (\ref{LAM1}) is algebraic and
the Galois group is finite.
Case $(i)$ splits in the following two subcases \cite{POO,WW}.
\begin{itemize}

\item[$(i.1)$]
The Lam\'e case. For a fixed integer $n$, this equation admits a solution (called \emph{Lam\'e function}) of the form
\begin{equation}
E(x)=\prod_{i=1}^3(x-e_i)^{k_i} P_m(x), \label{L1}
\end{equation}
being $P_m$ a monic polynomial of degree $m=n/2-(k_1+k_2+k_3)$ and $k_i \in \left\{0, \frac{1}{2} \right\}$, $i=1,2,3$. Since $m \in \N$,
eight different possibilities regarding $n$ appear: If $n$ is even, we have $k_1=k_2=k_3=0$ or just one zero $k_i$; if $n$ is odd, we could have all non-zero $k_i$'s or combinations with exactly one non-zero $k_i$. All these possibilities give rise to classes of Lam\'e functions.
Concerning parameter $B$, it must be one of the $m+1$ different roots $B_1,\ldots, B_{m+1}$
of certain irreducible polynomial of degree $m+1$, with all its roots real and simple~\cite{HAL}.
Furthermore, the numbers $B_i$  are reals.

\item[$(i.2)$] The Hermite case. Here we are not in case (i.1) and $n$ is an arbitrary natural number. We are also fall  in case 1 of Kovacic algorithm, but with a diagonal Galois group.

\end{itemize}

\begin{remark} \label{remlamesol} We notice that the polynomial $P_m$ in (i.1) satisfies a second order linear differential equation similar to the one that appears in the first case of Kovacic algorithm. In fact it is possible to obtain the above passing into normal form and applying Kovacic algorithm. Therefore the second linear independent solution is not algebraic and the associated Riccati equation has no rational first integral.
\end{remark}

\section*{Acknowledgements}
Remark \ref{remarkfactor} comes from a question posed by J.
Llibre to the third author some years ago. The authors are indebted
to J.A. Weil for several discussion and suggestions. The last
author also thanks J. Llibre  for many interesting discussions on
several subjects of this work.

\end{document}